\documentclass[10pt,a4paper]{article}
\makeatletter
\usepackage[dvipsnames]{xcolor}
\usepackage{lipsum}
\usepackage{tipa}
\usepackage{tikz}
\usepackage{tikz-cd}
\usetikzlibrary{matrix,arrows,decorations.pathmorphing}
\newcounter{magicrownumbers}

\usepackage[left=3cm,top=3cm,bottom=3cm,right=3cm,head=1cm,foot=1cm]{geometry}
\usepackage[pdftex,bookmarks=true,bookmarksnumbered=true]{hyperref}
\hypersetup{linkbordercolor=ForestGreen!25, citebordercolor=Goldenrod!60}
\usepackage{graphicx}
\usepackage{stmaryrd}
\usepackage{twoopt}
\usepackage{amssymb}
\usepackage{amsmath}
\usepackage{amsthm}
\usepackage{mathabx}
\usepackage{amsfonts}
\usepackage{amscd}
\usepackage{url}
\usepackage[all]{xy}
\newcommand{\ch}{\mathrm{ch}}

\newcommand{\tr}{\mathrm{tr}}

\newcommand{\hil}{\mathcal{H}}

\newcommand{\cl}{\mathrm{cl}}

\newcommand{\5}{\hspace{0,5cm}}
\newcommand{\3}{\vspace{0,3cm}}

\title{Analysis of Twisted Supercharge Families on Product Manifolds}
\date{31.3.2015}
\author{Antti J. Harju}
\begin{document}

\maketitle
\begin{abstract}
Twisted supercharge families on product manifolds $\mathbb{T} \times M$ have been applied in the analysis of the odd twisted K-theory. We shall suspend these families to the even twisted K-theory and solve their twisted families index problem. This is applied to give analytic representatives of the twisted K-theory classes on tori - including all the torsion classes. \3

\end{abstract}

\section*{Introduction}

Methods in Hamiltonian quantization were applied in the reference \cite{HM12} to give explicit constructions of gerbes and elements in the odd twisted K-theory groups of compact product manifolds $\mathbb{T} \times M$. The decomposable cup product twistings were considered. These 3-cohomology twisting classes are represented by $\tau = \alpha \smile \beta$ where $\alpha$ represents the generator of $H^1(\mathbb{T}, \mathbb{Z})$ and $\beta$ is a nontrivial cocycle in $H^2(M, \mathbb{Z})$. The twisted K-theory groups on the compact product manifolds are isomorphic to the group extensions 
\begin{eqnarray}\label{twistedk}
K^*(\mathbb{T} \times M, \tau) \simeq \{x \in K^*(M) : x = \lambda \otimes x\} \oplus_{\mu} \frac{K^{*-1}(M)}{(1-\lambda) \otimes K^{*-1}(M)}
\end{eqnarray}
for some extension cocycle $\mu$ where $\lambda$ denotes the K-theory class of a complex line bundle which is characterized by $\beta$ in cohomology. The methods of \cite{HM12} are relevant for the study of the second summand in \eqref{twistedk}. In this formalism, the operators which are associated to the twisted K-theory elements are called supercharges. The construction is essentially based on the theory of projective representations of the loop algebra of the circle group, $\mathfrak{lt} = \mathfrak{lu(1)}$, on Fock spaces.

Projective families of Dirac operators have been used for analytic realizations of the even twisted K-theory elements, see \cite{MMS05}, \cite{MMS09}, \cite{BG11}. Using the methods of \cite{MMS09} and \cite{BG11}, one can realize the even twisted K-theory classes with decomposable twisting in terms of projective families of pseudodifferential operators. Especially, the twisted characteristic classes for the twisted index bundles have been found. The characteristic classes are elements in the twisted cohomology theory which is computed from a de Rham complex with a differential modified to be compatible with the twistings. 

In \cite{HM12} we approached the problem to define the cohomology differently. In this approach, one pulls a gerbe and a supercharge to the covering space $\mathbb{R} \times M$ of $\mathbb{T} \times M$ where the gerbe trivializes and computes a characteristic form for the usual K-theoretic index problem (without any twist). These forms are not periodic in $\mathbb{R}$ and thus cannot be pulled back to the base $\mathbb{T} \times M$. However, one can project them to a quotient group of the odd cohomology group of the base: the cohomology theory needs to be made compatible  with the second summand of \eqref{twistedk}. We observed in \cite{HM12} that this construction is sensitive to some torsion information in the twisted K-theory. 

The goals of this work are to develop the supercharge approach to the twisted K-theory further and to connect it to the Dirac approach. In general, the odd twisted K-theory of a space $X$ can be seen as an image of a suspension map in the even twisted K-theory of $\mathbb{T} \times X$. This way we can transfer the supercharges to the even twisted K-theory. As an application of the representation theory of the loop algebra of the two-torus, we will realize these classes as supercharges based on the representations of $\mathfrak{lt}^2$, the Lie algebra of a loop group of the 2-torus. The methods to develop characteristic classes in \cite{HM12} will be extended to the even theory and the characteristic classes of the indexes of the suspended supercharges will be computed. In addition, by applying the methods of \cite{BG11}, both odd and even index characters for the supercharges will be computed in the twisted cohomology theory as well. 

The constructions in this work are based on the theory unbounded Fredholm operators. The twisted K-theory classes are associated to them explicitly by mapping them to the space of bounded operators under the approximate sign map. The topology in the space of unbounded Fredholm operators is induced by this map. The essential definitions in  the twisted K-theory and twisted cohomology are collected in the Appendix. We start with a short review of the superconnection methods in 1.1 - 1.4 which makes this work self contained. In 1.5 we formulate a superconnection formalism for the odd twisted K-theory in the twisted cohomology. The analogous construction in the even twisted K-theory was done in \cite{BG11}. We end the section 1 by solving an index character for the odd supercharges in the twisted cohomology. In section 2 and 3 the supercharge techniques are applied in the even twisted K-theory. The index characters for the suspended supercharges will be analyzed. In Section 4 we apply the techniques developed to give analytic constructions for the even and the odd twisted K-theory elements in the second summand of \eqref{twistedk} in the case of tori. \3

\noindent \textbf{Notation.} We shall use the twisted K-theory in the local formulation and so we need to fix an open cover for $\mathbb{T} \times M$ where $\mathbb{T}$ denotes the unit circle. Let $\{ \mathbb{T}_{\uparrow}, \mathbb{T}_{\downarrow}\}$ denote an open cover of $\mathbb{T}$ such that $\mathbb{T}_{\uparrow} \cap \mathbb{T}_{\downarrow}$ is a union of two disconnected components $(\mathbb{T}_{\uparrow} \cap \mathbb{T}_{\downarrow})^{(\pm 1)}$ containing $\pm 1 \in \mathbb{T}$. One views $\mathbb{T}_{\uparrow}$ as an upper hemisphere and $\mathbb{T}_{\downarrow}$ as a lower hemisphere of $\mathbb{T}$. Consider the open cover $\{\mathbb{T}_{\uparrow} \times M,\mathbb{T}_{\downarrow} \times M\}$ for $\mathbb{T} \times M$.  We shall often consider some analytic data (smooth functions, differential forms, operator families) which are defined locally over $\mathbb{T}_{\uparrow} \times M$ and then we use the symbol $A_{\uparrow}$ to indicate that the domain of $A_{\uparrow}$ is $\mathbb{T}_{\uparrow} \times M$, and similarly for $A_{\downarrow}$. The symbol $A_{\updownarrows}$ means that both components $A_{\uparrow}$ and $A_{\downarrow}$ are considered. We fix an open good cover for $M$. This cover will be denoted by $\{U_a: a \in I\}$. Since $M$ is compact, the index set $I$ can be taken to be finite. We shall not write the index set explicitly. 

We shall use the symbol $\ch(F_{\xi})$ to denote the Chern character form $\tr \exp(-\frac{F_{\xi}}{2 \pi i})$ if $F_{\xi}$ is a curvature of the complex vector bundle $\xi$.  We also use the normalization $\varphi$ for the differential forms on a manifold $X$ so that $\varphi$ sends the even form $\Omega \in \Lambda^{2n}(X)$ to $(2 \pi i)^{-n/2} \Omega$ and the odd form $\Phi \in \Lambda^{2n+1}(X)$ to $(2 \pi i)^{-n/2} \Phi$. \3
 
\noindent \textbf{Acknowledgment.} I wish to thank Danny Stevenson for help, especially for drawing my attention to the reference \cite{AS04}. 

\section{Character Maps for $K^1(\mathbb{T} \times M, \tau)$}

\noindent \textbf{1.1.}  Let $M$ denote a compact smooth manifold. A decomposable 3-cohomology class $\tau = \alpha \smile \beta$ of $\mathbb{T} \times M$ trivializes when restricted to the components of the cover $\{\mathbb{T}_{\uparrow} \times M, \mathbb{T}_{\downarrow} \times M\}$. A gerbe associated with a decomposable Dixmier-Douady class can be viewed as a pair of locally defined Hilbert bundles with connections $\textbf{F}_{\uparrow}$ and $\textbf{F}_{\downarrow}$ over $\mathbb{T}_{\uparrow} \times M$ and $\mathbb{T}_{\downarrow} \times M$, which can be identified over $(\mathbb{T}_{\uparrow} \cap \mathbb{T}_{\downarrow})^{(-1)} \times M$ and are subject to the local isomorphisms 
\begin{eqnarray*}
u_a : \textbf{F}_{\downarrow} \otimes \lambda \rightarrow  \textbf{F}_{\uparrow}   \5 \text{on} \5 (\mathbb{T}_{\uparrow} \cap \mathbb{T}_{\downarrow})^{(1)} \times U_a
\end{eqnarray*}
where $\lambda$ denotes a nontrivial complex line bundle on $M$ so that $\beta$ is a representative of $\lambda$ in cohomology. These maps are only locally defined on $M$ since $\lambda$ is a nontrivial bundle. The curvature forms of the local Hilbert bundles satisfy
\begin{eqnarray*}
u_a^*(F_{\uparrow}) = F_{\downarrow} + F_{\lambda} 
\end{eqnarray*}
where $F_{\downuparrows}$ and $F_{\lambda}$ are the curvature forms in $\textbf{F}_{\downuparrows}$ and $\lambda$. There are local forms $\Omega_{\downuparrows}$ defined over $\mathbb{T}_{\downuparrows} \times M$ such that $\Omega_{\uparrow} = \Omega_{\downarrow}$ on $(\mathbb{T}_{\uparrow} \cap \mathbb{T}_{\downarrow})^{(-1)} \times M$ and
\begin{eqnarray*}
\Omega_{\downarrow} - \Omega_{\uparrow} =  F_{\lambda} \5 \text{on} \5 (\mathbb{T}_{\uparrow} \cap \mathbb{T}_{\downarrow})^{(1)} \times M. 
\end{eqnarray*}
The form defined by $d \Omega_{\downuparrows}/2 \pi i$ on $\mathbb{T}_{\downuparrows} \times M$ is now a global 3-form over $\mathbb{T} \times M$ which is the curvature of the gerbe $\{\textbf{F}_{\downuparrows}\}$ and we denote it by $H$. 

Concretely, we can take $\Omega_{\uparrow} = 0$ and $\Omega_{\downarrow} = f_{\downarrow} F_{\lambda}$ where $f_{\downarrow}$ is any smooth real valued function on $\mathbb{T}_{\downarrow} \times M$ which is constant over the component $M$ and depends on the component $\mathbb{T}_{\downarrow}$ according to 
\begin{eqnarray*}
f_{\downarrow} = \left\{ \begin{array}{ll}
1 & \5 \text{on} \5 (\mathbb{T}_{\uparrow} \cap \mathbb{T}_{\downarrow})^{(1)} \times M \\  
0 & \5 \text{on} \5 (\mathbb{T}_{\uparrow} \cap \mathbb{T}_{\downarrow})^{(-1)} \times M  
\end{array} \right.
\end{eqnarray*}
We also set $f_{\uparrow} = 0$ which is a function on $\mathbb{T}_{\uparrow} \times M$ and now  we can write $\Omega_{\downuparrows} = f_{\downuparrows} F_{\lambda}$. Then $df_{\downarrow}$ is a bump form in $\mathbb{T}$ with integral equal to $1$. Thus, $H$ is equal to $\frac{d \phi}{2 \pi} \wedge \frac{F_{\lambda}}{2 \pi i}$ in cohomology and the Dixmier-Douady class of the gerbe $\{\textbf{F}_{\downuparrows}\}$ is indeed a decomposable class.\3

\noindent \textbf{1.2.} An explicit description of the local Hilbert bundles of 1.1 in terms of bundles of Fock spaces was developed in \cite{HM12}. Let $\mathcal{H}$ be an infinite dimensional complex Hilbert space with a basis $\{u_i: i \in \mathbb{Z} \}$. Next we construct a Fock space $\mathcal{F}$. The vacuum vector $|0 \rangle$ of $\mathcal{F}$ is defined by
\begin{eqnarray*}
|0 \rangle = u_{-1} \wedge u_{-2} \wedge u_{-3} \wedge \cdots 
\end{eqnarray*}
The basis is given by  
\begin{eqnarray}\label{base}
u_{a_1} \wedge \cdots \wedge u_{a_k} \wedge |0 \rangle_{b_1, \ldots, b_l}, \5 a_1 > \ldots > a_k \geq  0 > b_1 > \ldots > b_l
\end{eqnarray}
where $|0 \rangle_{b_1, \ldots, b_k}$ is a vacuum vector with the states $u_{b_x}$ removed. The charge of the vector \eqref{base} is defined to be $k-l$. The finite charge states are dense in the Fock space and we write
\begin{eqnarray*}
\mathcal{F} = \widehat{\bigoplus_{k \in \mathbb{Z}}} \mathcal{F}^{(k)}.
\end{eqnarray*}
There are two operators of interest: the shift $S$ is a unitary operator which raises the charge of the states, 
\begin{eqnarray*}
S: \mathcal{F}^{(k)} \rightarrow \mathcal{F}^{(k+1)}
\end{eqnarray*}
and $N$ is an unbounded operator which computes the charge of the state, i.e. $N = k \textbf{1}$ in $\mathcal{F}^{(k)}$. More precise definitions are given in \cite{HM12} and a very careful treatment of the Fock space theory in \cite{Ott95}.  

Bundles of Fock spaces can be used to give an explicit geometric description for the gerbe discussed in 1.1. We will define $\textbf{F}_{\downuparrows}$ to be the local Fock bundles over $\mathbb{T}_{\downuparrows} \times M$ such that for each $k \in \mathbb{Z}$, their charge $k$ subbundles will be tensored by $\lambda^{\otimes k}$ on $M$. In this setup, a gerbe can be constructed by letting the unitary transformations $u_a^{-1}$ raise the charge of a Fock bundle by one and in this process one needs to create a copy of  $\lambda$ on $M$. To create a complex line bundle one needs to fix a point in each fibre. This can only be done locally since we take $\lambda$ to be nontrivial. Now fix a locally defined $\mathbb{T}$-valued section $s_a: U_a \rightarrow \lambda$. Then define the local unitary transformations
\begin{eqnarray}\label{locuni}
&&u_a^{-1}: (\mathbb{T}_{\uparrow} \cap \mathbb{T}_{\downarrow})^{(1)} \times U_a \rightarrow U(\mathcal{F}) \\ \nonumber
&&u_a^{-1}(\phi, p) = s_a(p)S
\end{eqnarray}
These local maps create a Fock state which has the topological type of $\lambda$ over $M$. The local bundles $\textbf{F}_{\downuparrows}$ are given the connections $\nabla_{\downuparrows} = N \nabla_{\lambda}$ where $\nabla_{\lambda}$ is a connection in $\lambda$. According to 1.1, the Dixmier-Douady class of the gerbe $\{\textbf{F}_{\uparrow}, \textbf{F}_{\downarrow}\}$ is the decomposable class $\frac{d \phi}{2 \pi} \wedge \frac{F_{\lambda}}{2 \pi i}$. 

We shall also apply an irreducible projective representation of the loop algebra $\mathfrak{lt}$ (the Lie algebra of the loop group $L \mathbb{T}$) on the Fock space $\mathcal{F}$. Under this representation, the generators of the loop algebra satisfy 
\begin{eqnarray*}
[e_n, e_m] = - n \delta_{n,-m}.
\end{eqnarray*} 
and the rules $S e_0 S^{-1} = e_0 - 1$ and $S e_n S^{-1} = e_n$ for all $n \neq 0$. The explicit formulas are given in the reference \cite{HM12}. The charge operator $N$ is represented on the Fock space by $e_0$. \3 

\noindent \textbf{1.3.} A supercharge can be understood as a Dirac operator on an infinite dimensional loop group $L\mathbb{T}$. Therefore, we apply the Clifford algebra $\text{cl}(\mathfrak{lt})$, which is generated by $\psi_n, n \in \mathbb{Z}$ and the generators are subject to the anticommutation relations $\{ \psi_n, \psi_m \} = 2 \delta_{n,-m}$. We consider a spinor module which is a complex irreducible vacuum representation for $\text{cl}(\mathfrak{lt})$. Let $\eta$ denote the vacuum vector. Then $\psi_n \eta = 0$ for all $n < 0$ and we make $\psi_0$ act on the vacuum as the identity. 

Let $\textbf{S}$ denote a globally trivial bundle of $\cl(\mathfrak{lt})$-spinor modules over $\mathbb{T} \times M$. Then we can form another gerbe given by the local Hilbert bunles $\{\textbf{S} \otimes \textbf{F}_{\downuparrows}\}$. The odd supercharge $Q$ is given by the local families of $\textbf{Fred}^{(1)}_{\Psi}$-operators
\begin{eqnarray*}
	&&Q_{\downuparrows} : \mathbb{T}_{\downuparrows} \times M \rightarrow \textbf{Fred}^{(1)}_{\Psi} \\ 
	&& Q_{\downuparrows} = \sum_k \psi_k \otimes e_{-k} +  \frac{\phi}{2 \pi} \psi_0 \otimes \textbf{1}.
\end{eqnarray*}
Here $\textbf{Fred}^{(1)}_{\Psi}$ is a space of unbounded self-adjoint Fredholm operators, see Appendix A.1.

We introduce another dependence on the topology of $M$ by tensoring the gerbe with an arbitrary finite rank vector bundle $\xi$ on $M$. Now $\{\textbf{S} \otimes \textbf{F}_{\downuparrows} \otimes \xi\}$ is another representative of the decomposable Dixmier-Douady class and the connections in these Hilbert bundles are given by $e_0 \nabla_{\lambda} + \nabla_{\xi}$ where $\nabla_{\xi}$ is a fixed connection in $\xi$. Now the supercharge acts as $Q_{\downuparrows} \otimes 1$.  However, we shall exploit the notation and not write the second tensor component explicitly. \3

\noindent \textbf{1.4.} In \cite{HM12} we constructed a characteristic map for the odd supercharge. This construction applies the standard superconnection analysis in the  classical K-theory. One first defines the covering map $\pi: \mathbb{R} \times M \rightarrow \mathbb{T} \times M$ by
\begin{eqnarray*}
\pi: (x,p) \mapsto ([x], p), \5 [x] = x \text{ mod } 2 \pi. 
\end{eqnarray*}
Then the supercharge and the gerbe will be lifted to the covering. Since the 1-cohomology of $\mathbb{R}$ is trivial, the gerbe $\{\pi^*(\textbf{S} \otimes \textbf{F}_{\downuparrows} \otimes \xi) \}$ has a trivial Dixmier-Douady class, and hence defines a Hilbert bundle. The lifted supercharge is given by 
\begin{eqnarray*}
\pi^*(Q) = \sum_k \psi_k \otimes e_{-k} +  \frac{x}{2 \pi} \psi_0 \otimes \textbf{1}, \5 x \in \mathbb{R}.
\end{eqnarray*}
 Note that we write the lifted supercharge without the subscripts $\downuparrows$ because the bundle on which it acts is now a global geometric object, a trivializable bundle of Hilbert spaces. Over the covering one picks a superconnection for the lifted supercharge $\pi^*(Q)$.  Associated to the supercurvature $\mathbb{F}$ of the superconnection one defines the index character $\textbf{ch-ind}_1(\mathbb{F})$ in the usual way as a supertrace of the super Chern character, see \cite{HM12}, \cite{BF86b}. We found that this form is not $2 \pi$-periodic in the $\mathbb{R}$ direction and therefore cannot be pulled back to the base $\mathbb{T} \times M$. More precisely, under the translations by $2 \pi$ in $\mathbb{R}$, the index  character gets multiplied by the Chern character of the twisting line bundle $\ch(F_{\lambda})$. A solution is to project the cohomology to the quotient 
\begin{eqnarray}\label{14}
\frac{\sqrt{\pi} H^{\text{odd}}(\mathbb{T} \times M, \mathbb{Q})}{\sqrt{\pi} \frac{d \phi}{2 \pi} \wedge (1 - \ch(F_{\lambda})) \wedge \ch(K^0(M))}.
\end{eqnarray}
Here we consider $ H^{\text{odd}}(\mathbb{T} \times M, \mathbb{Q})$ as an abelian group and do not give it a ring or $\mathbb{Q}$-module structure. The pullback of the index character by an arbitrary section $\psi$ of $\pi$ is then well defined. \3

\noindent \textbf{Theorem 1} (\cite{HM12}). The odd twisted index character $\underline{\textbf{ch-ind}}^{\tau}_1(\mathbb{F}) = \psi^*(\textbf{ch-ind}_1(\mathbb{F}))$ projected on the quotient \eqref{14} is independent on the choice of the section $\psi$ and it is represented by the form
\begin{eqnarray*}
\underline{\textbf{ch-ind}}^{\tau}_1(\mathbb{F}) = \sqrt{\pi} \frac{d \phi}{2 \pi} \wedge \ch(F_{\xi}). 
\end{eqnarray*}
The approximate sign of the supercharge $F_{\downuparrows}$ defines a class in $K^1(\mathbb{T} \times M, \tau)$ given by 
\begin{eqnarray*}
F_{\downuparrows} = \frac{Q_{\downuparrows}}{\sqrt{1+Q_{\downuparrows}^2}}.
\end{eqnarray*}

\noindent \textbf{1.5.} Here we construct an odd index character for the supercharge $Q$ in the twisted cohomology. An odd superconnection in a gerbe of the type 1.1 is a descent datum $\mathbb{A} = (\mathbb{A}_{\downuparrows}, \nabla_{\lambda}, \Omega_{\downuparrows})$ where $\mathbb{A}_{\downuparrows}$ are odd superconnections for the local bundles $\textbf{S} \otimes \textbf{F}_{\downuparrows} \otimes \xi$ which satisfy the relations
\begin{eqnarray}\label{15}
u_a^*(\mathbb{A}_{\uparrow}) = \mathbb{A}_{\downarrow} \otimes \text{id} + \text{id} \otimes \nabla_{\lambda} \5 \text{on} \5 (\mathbb{T}_{\downarrow} \cap \mathbb{T}_{\uparrow})^{(1)} \times U_a. 
\end{eqnarray}
The local superconnections $\mathbb{A}_{\downuparrows}$ are odd elements in $\Lambda^{*}(\mathbb{T}_{\downuparrows} \times M, \text{End}(\textbf{S} \otimes \textbf{F}_{\downuparrows} \otimes \xi))$ which satisfy the Leibnitz rule 
\begin{eqnarray*}
\mathbb{A}_{\downuparrows}(\Phi f) = (d  \Phi) \wedge f + (-1)^k \Phi \wedge \mathbb{A}_{\downuparrows} f
\end{eqnarray*}
for $\Phi \in \Lambda^k(\mathbb{T}_{\downuparrows} \times M)$ and $f \in \Lambda^*(\mathbb{T}_{\downuparrows} \times M, \textbf{S} \otimes  \textbf{F}_{\downuparrows}\otimes \xi)$. The space $\Lambda^*(\mathbb{T}_{\downuparrows} \times M, \text{End}(\textbf{S} \otimes \textbf{F}_{\downuparrows}\otimes \xi))$ is extended with a formal symbol $\chi$ with $\chi^2 = 1$ and the $\mathbb{Z}_2$ grading is set such that $\chi$ and the differential 1-forms are considered as odd degree elements whereas the endomorphisms of $\textbf{S} \otimes \textbf{F}_{\downuparrows}\otimes \xi$ are even. Moreover, $\chi$ anticommutes with the differential 1-forms and commutes with the even elements. 

The twisted supercurvature is the data $\mathbb{F}_{\downuparrows}$ which is defined over $\mathbb{T}_{\downuparrows} \times M$ by 
\begin{eqnarray*}
\mathbb{F}_{\downuparrows} = \mathbb{A}^2_{\downuparrows} + \Omega_{\downuparrows}. 
\end{eqnarray*}
Now $\mathbb{F}$ is a globally defined form on $\mathbb{T} \times M$. The twisted odd characteristic form of the supercurvature $\mathbb{F}$ is defined by 
\begin{eqnarray*}
\textbf{ch-ind}^{\tau}_1(\mathbb{F}) = \varphi(\text{sTr} (e^{- \mathbb{F}})).
\end{eqnarray*}
where sTr is the odd supertrace which applies the Hilbert space trace on the component linear in $\chi$.

Denote by $H^{\text{odd}}(\mathbb{T} \times M, H)$ the odd twisted cohomology group on $\mathbb{T} \times M$ associated with the gerbe curvature $H$ (see also Appendix A.2). \3

\noindent \textbf{Proposition 1.} The form $\textbf{ch-ind}^{\tau}_1(\mathbb{F})$ is a cocycle in the twisted cohomology. If $\mathbb{A}$ and $\mathbb{A}'$ are twisted superconnections with the supercurvatures $\mathbb{F}$ and $\mathbb{F}'$, and if the twisting classes are determined by the gerbe curvatures 
\begin{eqnarray*}
H = \frac{df_{\downuparrows} \wedge F_{\lambda}}{2 \pi i} \5 \text{and} \5 H' = \frac{df_{\downuparrows} \wedge F_{\lambda}'}{2 \pi i},
\end{eqnarray*}
then there is an isomorphism of twisted cohomology groups $H^{\text{odd}}(\mathbb{T} \times M, H) \rightarrow H^{\text{odd}}(\mathbb{T} \times M, H') $  which sends $\textbf{ch-ind}^{\tau}_1(\mathbb{F})$ to $\textbf{ch-ind}^{\tau}_1(\mathbb{F}')$. \3

\noindent \text{Proof.} Since $\mathbb{A}^2_{\downuparrows}$ are ordinary supercurvatures, and $\mathbb{A}^2_{\downuparrows}$ commutes with $\Omega_{\downuparrows}$, we get
\begin{eqnarray*}
(d+H) \varphi (\text{sTr}(e^{- \mathbb{F}})) = - \frac{d \Omega}{2 \pi i} \wedge \varphi(\text{sTr} ( e^{- \mathbb{F}})) + H \wedge \varphi( \text{sTr}(e^{- \mathbb{F}})) = 0. 
\end{eqnarray*}

Denote by $\mathbb{A}$ and $\mathbb{A}'$ a pair of superconnections. Similarly, $F_{\lambda}$ and $F_{\lambda}'$ are the curvatures of the twisting line bundles and $\nabla_{\lambda}$ and $\nabla_{\lambda}'$ are the connections. The curvatures of the gerbes are then 
\begin{eqnarray*}
H = \frac{df_{\downuparrows} \wedge F_{\lambda}}{2 \pi i} \5 \text{and} \5 H' = \frac{df_{\downuparrows} \wedge F_{\lambda}'}{2 \pi i} 
\end{eqnarray*}
where we use the notation of 1.1. To prove the invariance on the superconnection, we consider the manifold $\mathbb{R} \times \mathbb{T} \times M$ where we can connect two superconnections by a smooth path. Define a pullback gerbe $\{\text{pr}^* (\textbf{S} \otimes \mathbf{F}_{\downuparrows}) \}$ over $\mathbb{R} \times \mathbb{T} \times M$ where $\text{pr}: \mathbb{R} \times \mathbb{T} \times M \rightarrow \mathbb{T} \times M$ is the natural projection map. Consider a pair of superconnections $\mathbb{A}_{\downuparrows}$ and $\mathbb{A}'_{\downuparrows}$ and the following local forms
\begin{eqnarray*}
\overline{\mathbb{A}^2}_{\downuparrows}(t) = (1-t) \text{pr}^*(\mathbb{A}^2_{\downuparrows}) + t \text{pr}^* ((\mathbb{A}')^2_{\downuparrows}) - dt \wedge \text{pr}^*(\mathbb{A}_{\downuparrows} - \mathbb{A}'_{\downuparrows}).
\end{eqnarray*}
Then we set a globally defined form 
\begin{eqnarray*}
\overline{\mathbb{A}^2}_{\downuparrows}(t) + \overline{\Omega}_{\downuparrows}(t)
\end{eqnarray*}
where $\overline{\Omega}_{\downuparrows}$ is chosen by 
\begin{eqnarray*}
\overline{\Omega}_{\downuparrows}(t) = (1-t) \text{pr}^*(f_{\downuparrows} F_{\lambda}) +  t \text{pr}^* (f_{\downuparrows} F'_{\lambda}) - dt \wedge \text{pr}^*(f_{\downuparrows}  (\nabla_{\lambda} - \nabla'_{\lambda})) - t \text{pr}^* (df_{\downuparrows} \wedge (\nabla_{\lambda} - \nabla'_{\lambda})) 
\end{eqnarray*} 
and $f_{\downuparrows}$ are the functions defined in 1.1. The gerbe curvature on $\mathbb{R} \times \mathbb{T} \times M $ is the global form 
\begin{eqnarray*}
\frac{d \overline{\Omega}_{\downuparrows}}{2 \pi i} = \text{pr}^* \Big( \frac{df_{\downuparrows}  \wedge F_{\lambda}}{2 \pi i} \Big) = \text{pr}^*(H) .
\end{eqnarray*}
Since $\overline{\mathbb{A}^2}_{\downuparrows}(t)$ and $\overline{\Omega}_{\downuparrows}(t)$ commute, it follows that
\begin{eqnarray*}
(d + \text{pr}^*(H)) \varphi(\text{sTr} (e^{-\overline{\mathbb{A}^2} - \overline{\Omega}}) )= 0
\end{eqnarray*}
and so $\varphi( \text{sTr} (e^{-\overline{\mathbb{A}^2} - \overline{\Omega}}))$ is a twisted cohomology cocycle. Consequently,
\begin{eqnarray*}
d_{\mathbb{R}}(\varphi (\text{sTr} (e^{-\overline{\mathbb{A}^2} - \overline{\Omega}}))) &=& - (d_{\mathbb{T} \times M} + \text{pr}^*(H)) \varphi(\text{sTr} (e^{-\overline{\mathbb{A}^2} - \overline{\Omega}})).
\end{eqnarray*}
Integration over the interval $[0,1]$ in $\mathbb{R}$ gives
\begin{eqnarray*} 
 e^{\frac{ df_{\downuparrows} \wedge (\nabla_{\lambda} - \nabla'_{\lambda})}{2 \pi i}} \wedge \varphi(\text{sTr}(e^{-\mathbb{F}'})) - \varphi(\text{sTr}(e^{-\mathbb{F}})) &=& - (d_{\mathbb{T} \times M} + H) \int^1_0 i_{\frac{\partial}{\partial t}} \varphi (\text{sTr}(e^{-\overline{\mathbb{A}^2} - \overline{\Omega}})) dt.
\end{eqnarray*}
It follows that the wedge product operation sends the class of $\textbf{ch-ind}_1(\mathbb{F}')$ in the complex $(\Lambda^*(\mathbb{T} \times M), d+H')$ to the class of $\textbf{ch-ind}_1(\mathbb{F})$ in the complex $(\Lambda^*(\mathbb{T} \times M), d+H)$. Since we have
\begin{eqnarray*}
H = H' - d \Big( \frac{ df_{\downuparrows} \wedge (\nabla_{\lambda} - \nabla'_{\lambda})}{2 \pi i} \Big)
\end{eqnarray*}
this map is an isomorphism of differential complexes (see Appendix A.2)
\begin{eqnarray*}
(\Lambda^*(\mathbb{T} \times M), d + H') \rightarrow (\Lambda^*(\mathbb{T} \times M), d + H). 
\end{eqnarray*}
Consequently, under the identification of differential complexes associated with the curvatures $H'$ and $H$, the character is independent on the choice of superconnection. \5 $\square$ \3

\noindent \textbf{1.6.} Following the notation of 1.1 we fix the one parameter families of local  supercurvatures   
\begin{eqnarray*}
\mathbb{A}_{\downuparrows}(t)^2 + \Omega_{\downuparrows} = (\sqrt{t} \chi Q_{\downuparrows} + \nabla_{\downuparrows})^2 + \Omega_{\downuparrows}
\end{eqnarray*}
for $t > 0$ with $\Omega_{\downuparrows} = f_{\downuparrows} F_{\lambda}$. We also write $F_{\downuparrows} = (\nabla_{\downuparrows})^2 = e_0 F_{\lambda} + F_{\xi}$. Recall that $e_0 = N$ is the Fock space operator which computes the charge of the states. These are the curvature forms in $\textbf{F}_{\downuparrows}$ but $F_{\downuparrows}$ cannot be glued to a global form since the charge grows by $1$ under the translations around the circle. The squares of the local families $Q_{\downuparrows}$ can be written by \cite{HM12}
\begin{eqnarray*}
Q_{\downuparrows}^2 = l_0^s + l_0^e + (e_0 + \frac{\phi}{2 \pi})^2, \5 \phi \in \mathbb{T}_{\downuparrows} \\
\end{eqnarray*}
where we have used 
\begin{eqnarray}\label{l}
l_0^s =  \sum_{n \geq 1} n \psi_n \psi_{-n} \5 \text{and} \5  l_0^e = 2\sum_{k \geq 1} e_k e_{-k}.
\end{eqnarray}
The operator $l_0^s + l_0^e$ is a positive unbounded operator with a dense domain given by the finite particle states which can be viewed as a Laplacian of the supercharge $Q$. We define a pair of locally defined families of orthogonal projections $P_{\downuparrows}$ onto the subspaces  
\begin{eqnarray}\label{P}
\bigoplus_{k \in \mathbb{Z}} \eta \otimes S^k |0 \rangle \otimes \mathbb{C}^n
\end{eqnarray}
in the fibres of $\{\textbf{S} \otimes \textbf{F}_{\downuparrows} \otimes \xi\}$. These subspaces are the kernels of $l_0^s + l_0^e$. The vectors \eqref{P} are eigenspaces for the number operator $e_0$ with eigenvalues $\mathbb{Z}$ and each eigenvalue has multiplicity one: note that 
\begin{eqnarray*}
e_0( \eta \otimes S^k |0 \rangle \otimes v) = k (\eta \otimes S^k |0 \rangle \otimes v ) 
\end{eqnarray*}
on the fibres. The symbol $\delta$ denotes the Dirac delta distribution. \3

\noindent \textbf{Proposition 2.} The form $\text{sTr}(e^{- \mathbb{A}^2_{\downuparrows}(t) - \Omega_{\downuparrows}})$ has an $\infty$-time limit given by the distribution valued odd differential form 
\begin{eqnarray}\label{delta}
 \lim_{t \rightarrow \infty} \text{sTr}(e^{- \mathbb{A}^2_{\downuparrows}(t) - \Omega_{\downuparrows}}) = \sqrt{\pi} \text{Tr}\big[ P_{\downuparrows} \delta( e_0 + \frac{\phi}{2 \pi}) \frac{d \phi}{2\pi} \wedge e^{-(F_{\xi} + e_0 F_{\lambda} + \Omega_{\downuparrows})}P_{\downuparrows}\big].
\end{eqnarray}
 
\noindent Proof. Since $\Omega_{\downuparrows}$ commutes with $(\sqrt{t} \chi Q_{\downuparrows} + \nabla_{\downuparrows})^2$ we can write 
\begin{eqnarray*}
\lim_{t \rightarrow \infty} \text{sTr}( e^{-\mathbb{A}_{\downuparrows}^2(t) - \Omega_{\downuparrows}}) = e^{ - \Omega_{\downuparrows}} \wedge \lim_{t \rightarrow \infty} \text{sTr}(e^{-\mathbb{A}^2_{\downuparrows}(t)}).
\end{eqnarray*}
The right hand side was computed in the reference \cite{HM12} which gives the proposition. \5 $\square$ \3

\noindent \textbf{Theorem 2.} If $Q_{\downuparrows}$ is a supercharge associated to the vacuum type $\xi$, and $\mathbb{F}$ is a global supercurvature form of $Q_{\downuparrows}$, then the character $\textbf{ch-ind}_1^{\tau}(\mathbb{F})$ is cohomologous to the twisted cohomology cocycle 
\begin{eqnarray}\label{ch}
\sqrt{\pi} \frac{d \phi}{2 \pi} \wedge \ch(F_{\xi}).
\end{eqnarray}

\noindent Proof. For the following analysis, recall that $\mathbb{T}_{\uparrow}$ and $\mathbb{T}_{\downarrow}$ denote the upper and lower hemisphere of $\mathbb{T}$ which we can parametrize using the angle variable $\phi$ so that $\mathbb{T}_{\uparrow} = (- \varepsilon, \pi + \varepsilon)$ and $\mathbb{T}_{\downarrow} = (\pi - \varepsilon, 2 \pi + \varepsilon)$ where $\varepsilon$ is a small real number. 

In the $\infty$-time limit, the form \eqref{delta} is nonzero on $\mathbb{T}_{\uparrow} \times M$ only if $e_0 + \frac{\phi}{2 \pi}$ has a zero eigenvalue in the image of $P_{\uparrow}$. This happens only at $\phi = 0$ and then $e_0$ is the zero operator on the eigenspace $\eta_0 \otimes |0 \rangle \otimes \mathbb{C}^n$. For $\phi = 0$ the form $-(F_{\xi} - \frac{\phi}{2 \pi} F_{\lambda} + \Omega_{\uparrow})$ is equal to $- F_{\xi}$. Over $\mathbb{T}_{\downarrow} \times M$, the eigenvalue of $e_0 + \frac{\phi}{2 \pi}$ is zero only if $\phi = 2 \pi$ because $e_0$ has the eigenvalue $-1$. The corresponding eigenspace is $\eta \otimes S^{-1} |0 \rangle \otimes \mathbb{C}^n$. For $\phi = 2 \pi$,  
\begin{eqnarray*}
-(F_{\xi} -  \frac{\phi}{2 \pi} F_{\lambda} + \Omega_{\downarrow}) = - F_{\xi} + F_{\lambda} - \Omega_{\downarrow} = - F_{\xi}
\end{eqnarray*}
Next we apply \eqref{delta}. Only the vectors $|0 \rangle$ and $S^{-1}|0 \rangle$ contribute to the traces which gives the global distribution valued form
\begin{eqnarray}\label{delta2}
\lim_{t \rightarrow \infty} \varphi (\text{sTr}(e^{-\mathbb{A}_{\downuparrows}^2 - \Omega_{\downuparrows}})) = \left\{ \begin{array}{ll} \sqrt{\pi} \delta(\frac{\phi}{2\pi}) \frac{d \phi}{2\pi} \wedge \ch(F_{\xi})  & \text{ if } \phi \in \mathbb{T}_{\uparrow}, \\
 \sqrt{\pi}  \delta(-1 + \frac{\phi}{2\pi}) \frac{d \phi}{2\pi} \wedge \ch(F_{\xi}) & \text{ if } \phi \in \mathbb{T}_{\downarrow}. \end{array} \right.
\end{eqnarray}
which localizes at $[\phi] = 0$. 

It follows from Proposition 1 that $\textbf{ch-ind}_1^{\tau}(\mathbb{F}(t))$ determines a constant path in the twisted cohomology group, and by \eqref{delta2}, the $t \rightarrow \infty$ limit is well defined as a distribution valued form. We need to understand the twisted cohomology class of $\textbf{ch-ind}_1^{\tau}(\mathbb{F}(t))$ for finite $t$ values. To this end, let us first consider the case where the twisting is turned off: $\lambda$ is chosen to be the trivial complex line bundle, $H$ and $\Omega_{\downuparrows}$ are zero forms. In this case, we can apply the character analysis as above which results in a path $\textbf{ch-ind}_1^0(\mathbb{F}(t))$ of odd differential forms. It determines a constant path in the ordinary odd cohomology of $\mathbb{T} \times M$. The limit \eqref{delta2} is independent on the chosen decomposable twisting, and so the $t \rightarrow \infty$ limit of $\textbf{ch-ind}_1^0(\mathbb{F}(t))$ is given by \eqref{delta2}. We call a differential form $\Theta \in \Lambda^*(\mathbb{T} \times M)$ a decomposable form if $\Theta = f(\phi) d \phi \wedge G$ so that $f(\phi) d \phi \in \Lambda^1(\mathbb{T})$ and $G \in \Lambda^*(M)$. The form  $\textbf{ch-ind}_1^0(\mathbb{F}(t))$ is a decomposable form for all $t > 0$: by to the subsection 6.2 in the reference \cite{HM12} we have
\begin{eqnarray*}
\text{sTr}( e^{-\mathbb{A}_{\downuparrows}^2(t)}) &=& \text{sTr}\Big( e^{-tQ_{\downuparrows}} + \sqrt{t} \sum_{n \geq 1} \frac{\chi \psi_0 d \phi \wedge (-F_{\downuparrows})^{n-1}}{2 \pi (n-1)!}  e^{-tQ_{\downuparrows}} + \sum_{n \geq 1} \frac{(- F_{\downuparrows})^n}{n!} e^{-tQ_{\downuparrows}} \Big) \\
&=& \frac{d \phi}{2\pi} \wedge \text{Tr} \Big( \sqrt{t} \sum_{n \geq 1} \frac{\psi_0 (-F_{\downuparrows})^{n-1}}{(n-1)!}  e^{-tQ_{\downuparrows}} \Big)
\end{eqnarray*}
because the terms that do not have a $\chi$ dependence are in the kernel of the supertrace. If $\Theta$ is a decomposable form and if $H$ is a decomposable twisting, then $H \wedge \Theta = 0$ because $H \wedge d \phi = 0$. So, if $\Theta$ is a closed decomposable form, then $\Theta$ defines a class in the twisted cohomology complex associated with $H$. The rest of the proof will rely on the following Lemma. 

\3 \noindent \textbf{Lemma 1.} If $\Theta$ and $\Theta'$ are closed decomposable forms of degree $2k +1$, and if the classes of $\Theta$ and $\Theta'$ in the ordinary cohomology complex are equal, then the classes of $\Theta$ and $\Theta'$ in the twisted cohomology complex are equal. 

\3 \noindent Proof. We have $\Theta - \Theta' = d \Phi$ for some $\Phi \in \Lambda^{2k}(\mathbb{T} \times M)$. Observe that $\Phi$ can be written as
\begin{eqnarray*}
\Phi = \sum_{k \in \mathbb{Z}} c_k e^{i k \phi} F_k + d \phi \wedge \Phi'
\end{eqnarray*} 
where $c_k$ are scalars, $c_0 = 0$, $\Phi' \in \Lambda^{2k-1}(\mathbb{T} \times M)$ and $F_k$ are  closed forms in $\Lambda^{2k}(M)$. The bump form $df_{\downuparrows}$ is the generator of the cohomology group $H^1(\mathbb{T})$. It follows that $e^{ik \phi} df_{\downuparrows}$ is an exact form for all $k \in \mathbb{Z} - \{0\}$ and therefore 
\begin{eqnarray}\label{G}
H \wedge \Phi = H \wedge \sum_{k \in \mathbb{Z}} c_k e^{i k \phi} F_k = d (\sum_{k \in \mathbb{Z}} g_k F_k \wedge \frac{F_{\lambda}}{2 \pi i}) := d G^{(1)}
\end{eqnarray}
for some smooth functions $g_k$ on $\mathbb{T}$, and therefore 
\begin{eqnarray*}
d \Phi = (d + H) \Phi - H \wedge \Phi = (d + H)(\Phi - G^{(1)}) + H \wedge G^{(1)}.
\end{eqnarray*}
The form $H \wedge G^{(1)}$ has degree $2k + 5$. From \eqref{G} we see that $G^{(1)}$ is the form 
\begin{eqnarray*}
G^{(1)} = \sum_{k \in \mathbb{Z}} c'_k e^{i k \phi} F'_k 
\end{eqnarray*}
for some scalars $c'_k$ with $c'_0 = 0$, and closed forms $F'_k \in \Lambda^{2k+2}(M)$. The exactness of $e^{ik \phi} df_{\downuparrows}$ for $k \in\mathbb{Z} - \{0\}$ implies that
\begin{eqnarray*}
H \wedge G^{(1)}  = d (\sum_{k \in \mathbb{Z}} g'_k F'_k \wedge \frac{F_{\lambda}}{2 \pi i}) := d G^{(2)}
\end{eqnarray*}
for some smooth functions $g'_k$ on $\mathbb{T}$. Now we can write 
\begin{eqnarray*}
d \Phi = (d + H)(\Phi - G^{(1)}) + H \wedge G^{(1)} = (d + H)(\Phi - G^{(1)} + G^{(2)}) - H \wedge G^{(2)}.
\end{eqnarray*}
and $H \wedge G^{(2)} \in \Lambda^{2k + 7}(\mathbb{T} \times M)$. Continuing inductively we find that for any $l \in \mathbb{N}$, there are the forms $G^{(i)} \in \Lambda^{2(k + i)}$ for $1 \leq i \leq l$ so that
\begin{eqnarray*}
d \Phi =  (d + H) \Big[ \Phi + \sum_{i = 1}^l (-1)^i G^{(i)} \Big] - (-1)^l H \wedge G^{(l)}.
\end{eqnarray*}
The degree of $H \wedge G^{(l)}$ exceeds the dimensionality of $\mathbb{T} \times M$ eventually. Thus, $d \Phi = (d + H) \Omega$ for some $\Omega \in \Lambda^{\text{even}}(\mathbb{T} \times M)$. \5 $\square$

\3 \noindent Proof of Theorem 2 continues. Let us first make the following hypothesis: for finite $t$-values, $\textbf{ch-ind}_1^{0}(\mathbb{F}(t))$ is represented by $\sqrt{\pi} \frac{d \phi}{2 \pi} \wedge \ch(F_{\xi})$ in the ordinary odd cohomology group. If the decomposable closed 3-form $H$ is fixed, then we have seen that $\textbf{ch-ind}_1^{0}(\mathbb{F}(t))$ determines a class in the odd twisted cohomology group associated with $H$ for any $t > 0$. Under our hypothesis, it follows from Lemma 1 that \eqref{ch} is a representative of the twisted cohomology class of $\textbf{ch-ind}_1^{0}(\mathbb{F}(t))$. On the other hand, if $\textbf{ch-ind}_1^{\tau}(\mathbb{F}(t))$ denotes the character associated with the $H$-twisting, then $\textbf{ch-ind}_1^{0}(\mathbb{F}(t))$ and $\textbf{ch-ind}_1^{\tau}(\mathbb{F}(t))$ determine the same twisted cohomology class for finite $t$-values because their $t \rightarrow \infty$ limits are equal, and in both cases, the twisted cohomology class is $t$-independent. Therefore \eqref{ch} is a representative of $\textbf{ch-ind}_1^{\tau}(\mathbb{F}(t))$. So, it is sufficient to check that $\textbf{ch-ind}_1^{0}(\mathbb{F}(t))$ is represented by $\sqrt{\pi} \frac{d \phi}{2 \pi} \wedge \ch(F_{\xi})$ in the ordinary cohomology for finite $t$-values. 

K{\" u}nneth's formula gives the isomorphism 
\begin{eqnarray*}
H^{\text{odd}}(\mathbb{T} \times M) = (H^0(\mathbb{T}) \otimes H^{\text{odd}}(M)) \oplus (H^1(\mathbb{T}) \otimes H^{\text{even}}(M)).
\end{eqnarray*} 
We can consider the closed decomposable forms of odd degree as representatives of the elements in the subgroup $H^1(\mathbb{T}) \otimes H^{\text{even}}(M)$. This subgroup is isomorphic to $H_0^{\text{cur}}(\mathbb{T}) \otimes H^{\text{even}}(M)$ where $H_0^{\text{cur}}(\mathbb{T})$ is the de Rham homology group of currents of degree $1$ on $\mathbb{T}$. Concretely, this isomorphism is the identity map on the $H^{\text{even}}(M)$ component and maps a closed de Rham 1-form $\Theta$ on $\mathbb{T}$ to the current $\hat{\Theta}$ defined by
\begin{eqnarray*}
\hat{\Theta}(\alpha) = \int_{\mathbb{T}} \alpha \Theta
\end{eqnarray*}
for all degree zero forms $\alpha \in C^{\infty}(\mathbb{T})$. Now $\textbf{ch-ind}_1^{0}(\mathbb{F}(t))$ determines a class in $H_0^{\text{cur}}(\mathbb{T}) \otimes H^{\text{even}}(M)$ which is independent of $t$. The $t \rightarrow \infty$ limit is well defined in $H_0^{\text{cur}}(\mathbb{T}) \otimes H^{\text{even}}(M)$ and it is equivalent to $\sqrt{\pi} \hat{d \phi}/2 \pi \otimes \ch(F_{\xi})$. So, at finite times, the cohomology class of $\textbf{ch-ind}_1^{0}(\mathbb{F}(t))$ is represented by $\sqrt{\pi} \frac{d \phi}{2 \pi} \wedge \ch(F_{\xi})$. \5 $\square$\3

The character $\textbf{ch-ind}_1^{\tau}(\mathbb{F})$ depends only on the twisted K-theory class of the supercharge $Q$. Namely, if we have a fixed gerbe and two supercharges which can be connected by a smooth path, then the superconnections can be connected by a smooth path, and as seen in the proof of Proposition 1, they define same class in the twisted cohomology.\3

\noindent \textbf{1.7.} The character forms introduced in 1.4 and 1.6 are significantly different since the two twisted cohomology theories applied cannot be identified in general. Let us consider an extreme case where they can be identified. Suppose that the cohomology of $M$ has a nontrivial torsion component in the degree 2-cohomology and $\lambda$ is a nontrivial flat line bundle associated with a torsion class. We choose a connection in $\lambda$ so that its curvature is the zero form. Then $\ch(F_{\lambda}) = 1$ and consequently, the normal subgroup $(1 - \ch(F_{\lambda})) \wedge \ch(K^0(M))$ in 1.4 is trivial and the twisted differential $d + H$ is just the ordinary de Rham differential $d$. Now the character in both approaches gets values in the usual de Rham complex $\Lambda^{\text{odd}}(\mathbb{T} \times M)$ and the cohomology class coincides with    
\begin{eqnarray*}
\sqrt{\pi} \frac{d \phi}{2 \pi} \wedge \ch(F_{\xi}). 
\end{eqnarray*}
This is an odd version of the local Mathai-Melrose-Singer twisted index theorem on a product with a torsion Dixmier-Douady class \cite{MMS05}. 

Now consider another extreme: the torsion part of the K-theory group $K^0(M)$ is trivial. However, given a nontrivial element $\lambda$ in $K^0(M)$, the summand
\begin{eqnarray*}
\frac{K^0(M)}{(1 - \lambda) \otimes K^0(M)}
\end{eqnarray*}
in  the twisted K-theory has a torsion component. The twisted cohomology theory cannot have torsion subgroups and the torsion information in the twisted K-theory will be lost when mapped to the twisted cohomology. However, when the character is evaluated in the quotient \eqref{14} all the torsion will be preserved. In fact, since the Chern character $\ch$ has a trivial kernel if $K^0(M)$ is torsion free, and since the character $\underline{\textbf{ch-ind}}_1^{\tau}(\mathbb{F})$ depends on $K^0(M)$ exactly through the Chern character, we see that the character $\underline{\textbf{ch-ind}}_1^{\tau}(\mathbb{F})$ distinguishes all the different twisted K-theory classes.

\section{Suspension in Twisted K-Theory}

\noindent \textbf{2.1.} We will follow the conventions of \cite{APS73} and define the suspension homotopy equivalence $\textbf{Fred}^{(1)} \rightarrow \Omega \textbf{Fred}^{(0)}$ by
\begin{eqnarray*}
	\text{susp}(A) = \left\{ \begin{array}{ll} \cos(s) + i \sin(s) A & \text{ if } s \in [0, \pi] \\
	   \cos(s) + i\sin(s) & \text{ if } s \in [\pi, 2 \pi]. \end{array} \right.
\end{eqnarray*}
When applied fiberwise in a bundle of Fredholm operators, one gets a homomorphism 
\begin{eqnarray*}
\text{susp}: K^1(\mathbb{T}_{\phi} \times M, \tau) \rightarrow K^0(\mathbb{T}_s \times \mathbb{T}_{\phi} \times M, \tau)
\end{eqnarray*}
with a trivial kernel. The Dixmier-Douady class of the gerbe on $\mathbb{T}^2 \times M$, which is still denoted by $\tau$, is the pullback of the Dixmier-Douady class in $\mathbb{T}_{\phi} \times M$ through the obvious projection $\mathbb{T}^2 \times M \rightarrow \mathbb{T}_{\phi} \times M$.\3

\noindent \textbf{2.2.} Let $Q_{\downuparrows}$ be the odd supercharge defined in 1.3. In order to study the index problem of the suspension of $Q_{\downuparrows}$, we expand the Fock spaces by tensoring them with an irreducible representation for the complexification of the two dimensional Clifford algebra generated by $\gamma^0,\gamma^1$ and subject to $\{\gamma^i, \gamma^j\} =  2\delta^{ij}$. An irreducible representation is unique, up to isomorphisms, and can be chosen by 
\begin{eqnarray*}
\gamma^0 = \begin{pmatrix} 0 & 1 \\ 1 & 0 \end{pmatrix}, \5 \gamma^1 = \begin{pmatrix} 0 & -i \\ i & 0 \end{pmatrix}. 
\end{eqnarray*}
Then we get the local self-adjoint Fredholm families
\begin{eqnarray*}
\widetilde{Q}_{\downuparrows} =\begin{pmatrix} 0 & \widetilde{Q}_{\downuparrows}^- \\ \widetilde{Q}_{\downuparrows}^+ & 0 \end{pmatrix} = \left\{ \begin{array}{ll} \cos(s) \gamma^0 + \sin(s) \gamma^1 Q_{\downuparrows} & \text{ if } s \in [0, \pi] \\
  \cos(s) \gamma^0 + \sin(s)\gamma^1 & \text{ if } s \in [\pi, 2\pi]. \end{array} \right.
\end{eqnarray*}
Similarly we can construct the self-adjoint families $\widetilde{F}_{\downuparrows}$ from the approximate sign $F_{\downuparrows}$ of $Q_{\downuparrows}$.\3

\noindent \textbf{Proposition 3.} If $Q_{\downuparrows}$ is a supercharge on $\mathbb{T}_{\phi} \times M$, then
\begin{eqnarray*}
\text{susp}(F_{\downuparrows})  =  \left\{\begin{array}{ll} \text{cos}(s) + i \sin(s) F_{\downuparrows} & \text{ if } s \in [0, \pi] \\
  \text{cos}(s) + i \sin(s) & \text{ if } s \in [\pi, 2\pi]. \end{array} \right.
\end{eqnarray*}
defines a class in $K^0(\mathbb{T}^2 \times M, \tau)$.\3

\noindent \text{Proof.}  We consider $\text{susp}(F_{\downuparrows})$ as the lower-left off-diagonal block of
\begin{eqnarray*}
\widetilde{F}_{\downuparrows} = \left\{ \begin{array}{ll} \cos(s) \gamma^0  + \sin(s) \gamma^1  F_{\downuparrows} & \text{ if } s \in [0, \pi] \\
   \cos(s) \gamma^0 + \sin(s)  \gamma^1 & \text{ if } s \in [\pi, 2\pi]. \end{array}\right.
\end{eqnarray*}
It is sufficient to check the appropriate compactness and continuity properties for $\widetilde{F}_{\downuparrows}$, see Appendix A.1. We can use the properties of the approximate sign $F_{\downuparrows}$ of $Q_{\downuparrows}$ which is known to be well defined twisted K-theory element. For $s \in [0, \pi]$,
\begin{eqnarray*}
1 - \widetilde{F}_{\downuparrows}^2 = 1 - \cos^2(s) - \sin^2(s) F_{\downuparrows}^2  =  \sin^2(s)  (1 - F_{\downuparrows}^2).
\end{eqnarray*}
This is a compact operator since $1 - F_{\downuparrows}^2$ are compact local families on $\mathbb{T}_{\phi} \times M$. Moreover, it varies continuously on $\mathbb{T}^2 \times M$ because $1 - F_{\downuparrows}^2$ and $\sin^2(s)$ do. Therefore, the parametrix of $\widetilde{F}_{\downuparrows}$ can be chosen to be $\widetilde{F}_{\downuparrows}$ itself. The continuity of $\widetilde{F}_{\downuparrows}$ in the strong operator topology follows from the continuity of $F_{\downuparrows}$, $\sin(s)$ and $\cos(s)$. 

If $s \in [\pi, 2\pi]$, then $1 - \widetilde{F}_{\downuparrows}^2 = 0$ which is a compact and continuous family. \5 $\square$ \3 

\noindent \textbf{2.3.} Associated to the families $\widetilde{Q}_{\downuparrows}$, there is a twisted index character created as follows. Let $\pi: \mathbb{T}_s \times \mathbb{R} \times M \rightarrow \mathbb{T}_s \times \mathbb{T}_{\phi} \times M$ denote the covering map 
\begin{eqnarray*}
\pi: (s,x,p) \mapsto (s, [x], p), \5 [x] = x \text{ mod } 2\pi. 
\end{eqnarray*}
The Dixmier-Douady class trivializes on the cover and we denote by $\pi^*(\textbf{S} \otimes \textbf{F}_{\downuparrows} \otimes \xi)$ the lifted gerbe, which is a trivializable Hilbert bundle. Then we pull back $\widetilde{Q}_{\downuparrows}$ to the covering $\mathbb{T}_s \times \mathbb{R} \times M$ and compute the ordinary characteristic form of the index bundle of $\pi^*(\widetilde{Q}^+)$. As in 1.4, we write $\pi^*(\widetilde{Q})$ without the subscripts $\downuparrows$ because this is a family acting on a  Hilbert bundle. The character form is well defined on the cover since the gerbe trivializes, however, the form is not invariant under the $2 \pi$-translations in $\mathbb{R}$: the translation by $ 2 \pi$ raises the charge of each Fock state by one, and consequently the curvature $\pi^*(F) = \pi^*(F_{\xi} + e_0 F_{\lambda})$ of the lifted Fock bundle transforms according to   
\begin{eqnarray*}
\pi^*(F)(s,x,p) \mapsto \pi^*(F)(s,x + 2 \pi,p) - \pi^*(F_{\lambda})
\end{eqnarray*}
for all $(s,x,p) \in \mathbb{T}_s \times \mathbb{R} \times M$ (recall that $S e_0 S^{-1} =  e_0 -1$). The goal is to pull the character form of the index back to $\mathbb{T}^2 \times M$ by choosing a section for $\pi$. To make this possible we define the target of the pullback to be the quotient 
\begin{eqnarray}\label{evencoh}
\frac{H^{\text{even}}(\mathbb{T}^2 \times M, \mathbb{Q})}{\frac{d \phi}{2 \pi} \wedge (1- \ch(F_{\lambda})) \wedge \ch(K^1(\mathbb{T}_s \times M))}
\end{eqnarray}
We shall see that the pullback valued in this quotient is well defined and also independent on the choice of the section $\psi$. This approach is an analogue of the one developed in \cite{HM12}, and sketched in 1.4, in the odd twisted K-theory which now factorizes through the classical families index computation. 

 To compute an ordinary characteristic polynomial for the lifted index bundle of the Fredholm section $\pi^*(\widetilde{Q})$ we choose a one parameter family of Bismut superconnections
\begin{eqnarray*}
\mathbb{A}(t) = \sqrt{t} \pi^*(\widetilde{Q}) + \pi^*(\nabla).
\end{eqnarray*}
The supercurvature is given by
\begin{eqnarray}\label{supercurvature}
	\mathbb{F}(t) = t \pi^*(\widetilde{Q})^2 + \sqrt{t}[\pi^*(\nabla), \pi^*(\widetilde{Q})] + \pi^*(F),
\end{eqnarray}
The super Chern character of $\mathbb{A}(t)$ is defined by 
\begin{eqnarray}\label{charform}
	\textbf{ch-ind}_0(\mathbb{F}(t)) = \varphi(\text{sTr}(e^{- \mathbb{F}(t)})),
\end{eqnarray}
where sTr denotes the supertrace in the even formalism: it is defined to be $2^{-1}  \text{Tr}' \circ i\gamma^0 \gamma^1$ where $\text{Tr}'$ is the trace of the Hilbert space that is tensored by the finite Clifford module $\mathbb{C}^2$. We use $\text{Tr}$ to denote the Hilbert space trace of Section 1. The following transgression formula holds, \cite{BGV04} 
\begin{eqnarray*}
\frac{d}{dt} \text{sTr}(e^{- \mathbb{F}(t)})= - d( \text{sTr} \Big(\frac{d \mathbb{A}(t)}{dt} e^{- \mathbb{F}(t)}  \Big)).
\end{eqnarray*}
It follows that when pulled back to the base $\mathbb{T}^2 \times M$, the equivalence class of the character form in \eqref{evencoh} becomes homotopy invariant. 

We shall first compute the $\infty$-time limit of the character locally over $(0, \pi)$ in the $\mathbb{T}_s$ direction where the character localizes as a distribution valued form. Then we can extend it over $\mathbb{T}_s$. Notice that $\pi^*(\widetilde{Q}^2)$ is the identity operator for $t \in [\pi, 2 \pi]$, and so for these values of $s$, the heat operator $e^{- \mathbb{F}(t)}$ is not trace class. Denote by $P$ the family of projection operators on the cover $\mathbb{T}_s \times \mathbb{R} \times M$ onto the subspaces
\begin{eqnarray*}
\eta_0 \otimes \bigoplus_{k \in \mathbb{Z}} S^k |0 \rangle \otimes \mathbb{C}^n
\end{eqnarray*}
in the fibres of $\pi^*(\textbf{S} \otimes \textbf{F}_{\downuparrows} \otimes \xi)$.\3

\noindent \textbf{Proposition 4.} The $\infty$-time limit of $\textbf{ch-ind}_0(\mathbb{F}(t))$ exists as a distribution valued form on $(0, \pi) \times \mathbb{R} \times M$  and is equal to  
\begin{eqnarray*}
\text{Tr }( P\delta\Big( s - \frac{\pi}{2 } \Big) \delta\Big(e_0 + \frac{x}{2 \pi} \Big) \frac{ds \wedge dx}{2 \pi}  \wedge e^{-(2 \pi i)^{-1}(F_{\xi} + e_0 F_{\lambda})} P)
\end{eqnarray*}

\noindent Proof. Consider the supercurvature \eqref{supercurvature}. We can write $\pi^*(\nabla) = d_{s} + d_{x} + d_M + \Omega$ where $x$ is the coordinate of $\mathbb{R}$ and $\Omega$ is a matrix valued connection form. Then 
\begin{eqnarray*}
	[\pi^*(\nabla), \pi^*(\widetilde{Q})] = - \sin(s) \gamma^0 ds +  \cos(s) \gamma^1(\pi^*(\widetilde{Q}))ds + \sin(s) \gamma^1 \psi_0 \frac{dx}{2 \pi}. 
\end{eqnarray*}
In the following we shall also need
\begin{eqnarray*}
  [\pi^*(\nabla), \pi^*(\widetilde{Q})]^2 = -  \sin^2(s) \gamma^0 \gamma^1 \psi_0 ds \wedge \frac{dx}{\pi}. 
\end{eqnarray*}

Next we apply the Volterra series expansion to compute the exponential of the supercurvature $\mathbb{F}(t)$. Let $\mathbb{F}(t) = t((\pi^*(\widetilde{Q}))^2 + K_t)$ so that $K_t = \frac{1}{\sqrt{t}} [\pi^*(\nabla), \pi^*(\widetilde{Q})] + \frac{1}{t} \pi^*(F)$. Then
\begin{eqnarray*}
	e^{-\mathbb{F}(t)} &=& e^{- t \pi^*(\widetilde{Q})^2} \\ &+& \sum_{n \geq 1} (-t)^n \int_{\triangle_n} e^{-tx_1 \pi^*(\widetilde{Q})^2} K_t e^{-tx_2 \pi^*(\widetilde{Q})^2} \cdots e^{-tx_n \pi^*(\widetilde{Q})^2} K_t e^{-tx_{n+1} \pi^*(\widetilde{Q})^2} dx_1 \cdots dx_{n+1} \\
	&=& e^{-t\pi^*(\widetilde{Q})^2} + \sum_{n \geq 1} (-t)^n \int_{\triangle_n} K_t^n  e^{-t \pi^*(\widetilde{Q})^2} dx_1 \cdots dx_{n+1} \\
	&=& \sum_{n \geq 0} \frac{(-t)^n}{n!} K_t^n  e^{-t \pi^*(\widetilde{Q})^2},
\end{eqnarray*}
where the integration is over the standard $n$-simplex $\triangle_n$ in $\mathbb{R}^{n+1}$. We used the obvious identity $[\pi^*(\widetilde{Q}), K_t] = 0$. A direct computation gives           
\begin{eqnarray*}
	 K_t^n &=& \sum_{k = 0}^{n} \binom{n}{k} \Big(\frac{[\pi^*(\nabla), \pi^*(\widetilde{Q})]}{\sqrt{t}}  \Big)^{n-k} \wedge \Big(\frac{\pi^*(F)}{t}\Big)^{k} \\ 
	 &=& n(n-1) \frac{- \sin^2(s) \gamma^0 \gamma^1 \psi_0 ds \wedge dx }{\pi t} \wedge \Big( \frac{\pi^*(F)}{t}\Big)^{n-2} \\ 
	&+& n \frac{[\pi^*(\nabla), \pi^*(\widetilde{Q})]}{\sqrt{t}} \wedge \Big( \frac{\pi^*(F)}{t}\Big)^{n-1} + \Big(\frac{\pi^*(F)}{t}\Big)^{n}
\end{eqnarray*}
and it follows
\begin{eqnarray*}
	e^{-\mathbb{F}(t)} &=& t \sum_{n \geq 2}\frac{- \sin^2(s) \gamma^0 \gamma^1 \psi_0 ds \wedge dx  }{\pi} \wedge \frac{(- \pi^*(F))^{n-2}}{(n-2)!} e^{-t\pi^*(\widetilde{Q})^2} \\
	&-& \sqrt{t} \sum_{n \geq 1} [\pi^*(\nabla), \pi^*(\widetilde{Q})] \wedge \frac{(- \pi^*(F))^{n-1}}{(n-1)!} e^{-t\pi^*(\widetilde{Q})^2} +  \sum_{n \geq 0} \frac{(-\pi^*(F))^{n}}{n!} e^{-t\pi^*(\widetilde{Q})^2}.
\end{eqnarray*}
The second term in the sum is in the kernel of the supertrace since it is off diagonal with respect to the chiral grading. Also the third term vanishes because both of the diagonal components are equal. The square $\pi^*(\widetilde{Q})^2$ is equal to 
\begin{eqnarray*}
\pi^*(\widetilde{Q})^2 = \cos^2(s)  + \sin^2(s) \pi^*(Q)^2 = \cos^2(s) + \sin^2(s) (l_0^s + l_0^e + (e_0 + \frac{x}{2 \pi})^2),
\end{eqnarray*}
where $l_0^s$ and $l_0^e$ are given in \eqref{l} and the kernel of $l_0^s + l_0^f$ is the image of $P$ in the fibres of the Fock bundle. 

Up to contributions which vanish under the super trace, $e^{-\mathbb{F}(t)}$ is equal to the product of the endomorphism valued forms
\begin{eqnarray*}
&&\sqrt{\frac{t}{\pi}} e^{-t( \sin(s) (e_0 + \frac{x}{2 \pi}))^2}, \5 \sqrt{\frac{t}{\pi}} e^{-t \cos^2(s)} \5 \text{and} \\
	&& \sum_{n \geq 0} (- \sin^2(s) \gamma^0 \gamma^1 \psi_0) ds \wedge dx \wedge \frac{(- \pi^*(F))^{n}}{n!} e^{-t \sin^2(s)(l_0^s + l_0^e)}. 
\end{eqnarray*}
Since $s \in (0,\pi)$ it follows $\sin(s) \neq 0$. So, in the limit $t \rightarrow \infty$, only the kernel subspaces of $l_0^s + l_0^e$ contributes to the trace. Each of the endomorphism valued forms have $t \rightarrow \infty$ limit which are equal to 
\begin{eqnarray*}
	&& \delta\Big( \sin(s)  (e_0 + \frac{x}{2 \pi})  \Big), \5 \delta (\cos(s)) \5 \text{and} \\
	&&  P \sum_{n \geq 0} (- \sin^2(s) \gamma^0 \gamma^1 \psi_0 ) ds \wedge dx \wedge \frac{(- \pi^*(F))^{n}}{n!} P.
\end{eqnarray*}
Recall that $\psi_0$ is the identity operator on the image of $P$. Therefore, up to contributions which vanish under the supertrace
\begin{eqnarray*}
 \lim_{t \rightarrow \infty} \varphi (\text{sTr} (e^{-\mathbb{F}(t)})) \simeq - \text{sTr}(P \gamma^0 \gamma^1 \delta\Big( s - \frac{\pi}{2} \Big)   \delta\Big( e_0 + \frac{x}{2 \pi} \Big) \frac{ds \wedge dx}{2 \pi i} \wedge e^{-(2 \pi i)^{-1} \pi^* (F_{\xi} +  e_0  F_{\lambda})} P )
\end{eqnarray*}
where we have applied the normalization $\varphi$. When we substitute the even supertrace $\text{sTr} = 2^{-1} \text{Tr}' \circ i\gamma^0 \gamma^1$ we get the proposition. \5 $\square$ \3
 
The distribution valued form in the proposition localizes at $s = \pi/2$ and so we can extend it to  a distribution valued form on $\mathbb{T}^2 \times M$ by defining it to be equal to  zero for $s \in [\pi, 2 \pi]$. Since the kernel of the Fredholm section $\pi^*(\widetilde{Q})$ at these values of $s$ are the zero vector space, this extension is compatible with index theory. 

Fix a sections $\psi: \mathbb{T}^2 \times M \rightarrow \mathbb{T} \times \mathbb{R} \times M$ for the covering map $\pi$. Then define 
\begin{eqnarray*}
\underline{\textbf{ch-ind}}_0^{\tau}(\mathbb{F}) = \psi^*(\textbf{ch-ind}_0)(\mathbb{F}).
\end{eqnarray*}

\noindent \textbf{Theorem 3.} The index character $\underline{\textbf{ch-ind}}_0^{\tau}(\mathbb{F})$ defines a class in the quotient group \eqref{evencoh}. The class is independent on the choice of a superconnection and on the section $\psi$. If $Q_{\downuparrows}$ is a supercharge associated to a gerbe with a vacuum bundle $\xi$, then the index character is represented by
\begin{eqnarray}\label{form}
 \frac{d s}{2 \pi} \wedge \frac{d \phi}{2 \pi} \wedge \ch(F_{\xi}). 
\end{eqnarray}
If the even index characters of the suspended families of the odd supercharges $Q_{\downuparrows}$ and $Q'_{\downuparrows}$ are nonequivalent, then $\text{susp}(F_{\downuparrows}) \neq \text{susp}(F'_{\downuparrows})$ in $K^0(\mathbb{T}^2 \times M, \tau)$.\3

\noindent Proof. Consider the distribution valued form of Proposition 4. This form is nonzero only if $e_0 + \frac{x}{2 \pi}$ has a zero eigenvalue in the image of the projection $P$. The eigenvalues of $e_0$ are integers and so $x$ can be an integer times $2 \pi$. Then, at $x = 2 \pi k$ for $k \in \mathbb{Z}$, the kernel of $e_0 + \frac{x}{2 \pi}$ is given by $\eta_0 \otimes S^{-k} |0 \rangle \otimes \mathbb{C}^n$. Therefore, we can write 
\begin{eqnarray*}
\lim_{t \rightarrow \infty} \text{sTr}(e^{-\mathbb{F}(t)}) = \sum_{k \in \mathbb{Z}} \delta\Big( s - \frac{\pi}{2 } \Big) \delta\Big( k + \frac{x}{2 \pi} \Big) \frac{ds \wedge dx}{2 \pi}  \wedge \text{tr}_{\xi}( e^{-(2 \pi i)^{-1}(F_{\xi} - \frac{x}{2 \pi} F_{\lambda})})
\end{eqnarray*}
where $\text{tr}_{\xi}$ is the trace in the fibres $\mathbb{C}^n$ of $\xi$. This form is not invariant under translations by $2 \pi$ in the $\mathbb{R}$ direction since the curvature $F_{\xi} - \frac{x}{2 \pi} F_{\lambda}$ depends on $x$. However, the dependence on the Chern character of the twisting bundle $\lambda$ is eliminated in the quotient group \eqref{evencoh}. Therefore the pullback of the character defines a well defined element in the quotient group. For the same reason the choice of $\psi$ is arbitrary. The $\infty$-time limit of the character pulled to the base is given by
\begin{eqnarray*}
 \delta\Big( s - \frac{\pi}{2 } \Big) \delta\Big( \frac{\phi}{2 \pi} \Big ) \frac{ds \wedge d \phi}{2 \pi}  \wedge \ch(F_{\xi}). 
\end{eqnarray*}
The distribution forms $\delta\Big( s - \frac{\pi}{2 } \Big) ds$ and $\delta\Big( \frac{\phi}{2 \pi} \Big ) d \phi$ on the unit circles $\mathbb{T}_{s}$ and $\mathbb{T}_{\phi}$ have volumes equal to $1$ on $\mathbb{T}_s$ and $2 \pi$ on $\mathbb{T}_{\phi}$. Therefore we can use the argument of the proof of Theorem 2 to deduce that at finite times, these forms are equivalent to \eqref{form}. 

Suppose that the even index characters of $\widetilde{Q}_{\downuparrows}$ and $\widetilde{Q}'_{\downuparrows}$ are nonequal. It follows that there will be no homotopy connecting these local families, or equivalently, there will be no homotopy connecting the odd supercharges $Q_{\downuparrows}$ and $Q'_{\downuparrows}$. The topology in the space of unbounded operators is induced from the strong operator topology in $\mathcal{B}(\hil)$ and so the approximate signs $F_{\downuparrows}$ and $F'_{\downuparrows}$ cannot be homotopic either. Then $\text{susp}(F_{\downuparrows})\neq \text{susp}(F'_{\downuparrows})$ as elements in the even twisted K-theory. \5 $\square$ \3

\noindent \textbf{2.4.} Define a normalized desuspension map by
\begin{eqnarray*}
\Sigma^{-1} &=& \sqrt{\pi} \int_{\mathbb{T}_s} : \frac{H^{\text{even}}(\mathbb{T}^2 \times M, \mathbb{Q})}{\frac{d \phi}{2 \pi} \wedge (1- \ch(F_{\lambda})) \wedge \ch(K^1(\mathbb{T}_s \times M))} \\ &\rightarrow& \frac{\sqrt{\pi}H^{\text{odd}}(\mathbb{T}_{\phi} \times M, \mathbb{Q})}{\sqrt{\pi}\frac{d \phi}{2 \pi} \wedge (1- \ch(F_{\lambda})) \wedge \ch(K^0(M))}.
\end{eqnarray*}
Then we see that the odd index character of 1.4 factors through the even character of the twisted index bundle as
\begin{eqnarray*}
\underline{\textbf{ch-ind}}_1^{\tau} = \Sigma^{-1} \circ \underline{\textbf{ch-ind}}^{\tau}_0 \circ \tilde{\text{susp}}. 
\end{eqnarray*}
where $\tilde{\text{susp}}$ is the operation $Q_{\downuparrows} \mapsto \widetilde{Q}_{\downuparrows}$. The character for the odd index in the case of Dirac families is defined correspondingly, see \cite{BF86b}.

\section{Supercharge Suspension}

\noindent \textbf{3.1.} Next we apply the representation theory of $\mathfrak{lt}^2$ to construct another realization for the suspension $\text{susp}(F_{\downuparrows})$ in the twisted K-theory. 

Consider the complexified Clifford algebra, $\text{cl}(\mathfrak{lt}^2)$, of the real Clifford algebra of $\mathfrak{lt}^2$ which is the polynomial algebra generated by $\psi_n^i$ with $i = 0,1$ and $n \in \mathbb{Z}$ subject to the relations
\begin{eqnarray*}
	\{\psi_n^i, \psi_m^j \} = 2 \delta^{ij} \delta_{n,-m} . 
\end{eqnarray*}
We use a vacuum representation for $\text{cl}(\mathfrak{lt^2})$. The vacuum is a two dimensional subspace which is an irreducible module for the two dimensional subalgebra generated by $\psi_0^0$ and $\psi_0^1$. We choose the basis of the vacuum $\{ \eta_1, \eta_2 \}$ so that $\psi_0^0$ and $\psi_0^1$ restrict on the vacuum subspaces to the operators
\begin{eqnarray*}
	\psi_0^0 = \left(\begin{array}{cc}
		0 & 1 \\
		1 & 0 
	\end{array}\right) \5 \text{and} \5  \psi^1_0 = \left(\begin{array}{cc}
		0 & -i \\
		i & 0 
	\end{array}\right).
\end{eqnarray*}
The vacuum subspace is annihilated by the operators $\psi_n^i$ for all $n < 0$ with $i = 0,1$. The basis vectors can be created by acting on the vacuum by the subalgebra generated by $\psi_n^i$ with $n > 0$. Let us denote by $\mathcal{H}_s^e$ the spinor module obtained this way. 

Consider the operator 
\begin{eqnarray*}
f = \sum_{k > 0} (\psi^0_k \psi^0_{-k} + \psi^1_k \psi^1_{-k})
\end{eqnarray*}
which is a positive operator and diagonal in the basis of $\mathcal{H}_s^e$. The operator $f$ computes the number of the Clifford algebra generators needed to create the state from the vacuum. Then we define the self-adjoint operator with eigenvalues $\pm 1$ by 
\begin{eqnarray*}
\Gamma = - i \psi_0^0 \psi_0^1 (-1)^f. 
\end{eqnarray*}
The $+1$ eigenspace of $\Gamma$ in $\mathcal{H}_s^e$ is spanned by $\eta_1$, which is the $+1$ eigenspace of $- i \psi_0^0 \psi_0^1$, the states that are generated by acting with an even number of Clifford algebra elements on $\eta_1$ and the states that are generated by acting with an odd number of Clifford algebra elements on $\eta_2$. We define the polarization $\mathcal{H}^e_s = \mathcal{H}_s^+ \oplus \mathcal{H}_s^-$ with respect to the chiral grading determined by $\Gamma$.\3
 
\noindent \textbf{3.2.} Consider the unit circle $\mathbb{T}_s$ with an open cover given by $\{(\mathbb{T}_s)_{\downuparrows} \}$. Define the canonically trivial bundles of Fock spaces $\textbf{F}^s_{\downuparrows}$ on $(\mathbb{T}_s)_{\downuparrows}$ and introduce the constant family of Fock space shift operators 
\begin{eqnarray*}
&& g: ((\mathbb{T}_s)_{\uparrow} \cap (\mathbb{T}_s)_{\downarrow})^{(1)} \rightarrow U(\mathcal{F}); \\
&& g(s) =  S
\end{eqnarray*}
Then we glue $\textbf{F}^s_{\downuparrows}$ over $((\mathbb{T}_s)_{\uparrow} \cap (\mathbb{T}_s)_{\downarrow})^{(-1)}$ and use $g$ as a $U(\mathcal{F})$ valued transition function to glue the local bundles $\textbf{F}^s_{\downuparrows}$ over the domain of $g$. In particular, over the domain of $g$, a vector $\Psi$ in the fibre of $\textbf{F}^s_{\downarrow}$ is identified with the vector $g(s) \Psi = S \Psi$ in the fibre of $\textbf{F}^s_{\uparrow}$. This results in a bundle of Fock spaces over $\mathbb{T}_s$ which we denote by $\textbf{F}^s$. As explained in 1.3, now one can construct a family of supercharges on $\mathbb{T}_s$ with spectral flow around the circle. This would give the analytic realization for the generator of the ordinary K-theory group $K^1(\mathbb{T}_s)$. 

Let us then form the local tensor product bundles 
\begin{eqnarray*}
\textbf{F}^s \boxtimes \textbf{F}_{\downuparrows} \rightarrow \mathbb{T}_s \times (\mathbb{T}_{\phi})_{\downuparrows} \times M,
\end{eqnarray*}
whose fibres are copies of $\mathcal{F}^s \otimes \mathcal{F}$. The Dixmier-Douady class of $\{\textbf{F}^s \boxtimes \textbf{F}_{\downuparrows}\}$ is the product class $\tau$ pulled back to $\mathbb{T}^2 \times M$ through the projection $\pi$ onto $\mathbb{T}_{\phi} \times M$. The new bundle does not contribute to the Dixmier-Douady class because its transition functions act on a different tensor product component.

The construction above allows one to construct a new supercharge with spectral flows around both circles. In 2.3 we found that the suspended family is indeed associated with such an operator since the K-theoretic data was localized in the distribution forms in the direction of both unit circles.

Both tensor components $\mathcal{F}^s$ and $\mathcal{F}$ are Fock spaces with properties introduced in 1.2. On $\mathcal{F}^s \otimes \mathcal{F}$ we have two copies of projective representations for the loop algebra $\mathfrak{lt}$, given by $e_n \otimes \textbf{1}$ and $\textbf{1} \otimes e_n$ which we shall denote by $e_n$ and $f_n$. These operators define a projective representations of $\mathfrak{lt}^2$ and verify the commutation relations  
\begin{eqnarray*}
	[e_n, f_m ] = 0 ,\5 [e_n, e_m] = [f_n, f_m] = - n \delta_{n,-m}. 
\end{eqnarray*}

Let $\textbf{S}^e$ denote the globally trivial bundle of $\mathcal{H}_s^e$-spinor modules over $\mathbb{T}^2 \times M$. We realize the gerbe associated to the projective Dixmier-Douady class as the pair of Hilbert bundles
\begin{eqnarray}\label{evengerbe}
\{\textbf{S}^e \otimes (\textbf{F}^s \boxtimes \textbf{F}_{\downuparrows}) \otimes \xi \}
\end{eqnarray}
where $\xi$ is a finite rank complex vector bundle over $M$. \3

\noindent \textbf{3.3.} The even supercharge is a pair of locally defined self-adjoint unbounded Fredholm operator families 
\begin{eqnarray*}
	Q_{\downuparrows}^e: \mathbb{T}_s \times \mathbb{T}_{\downuparrows} \times M \rightarrow \textbf{Fred}^{(0)}_{\Psi}
\end{eqnarray*}
defined by
\begin{eqnarray*}
	Q_{\downuparrows}^e(s,\phi,p) = \sum_k \psi^0_k \otimes e_{-k} + \sum_k \psi^1_k \otimes f_{-k} + \frac{s}{2 \pi} \psi^0_0 \otimes \textbf{1} + \frac{\phi}{2 \pi} \psi^1_0  \otimes \textbf{1}. 
\end{eqnarray*}
where the coordinate $\phi$ gets values in $\mathbb{T}^{\uparrow}$ or $\mathbb{T}^{\downarrow}$ and $s \in \mathbb{T}_s$. 

The families $Q^{\downuparrows}$ are densely defined and their domain consists of the finite particle Fock space states. Moreover, they are antidiagonal with respect to the chiral grading so we can write 
\begin{eqnarray*}
Q_{\downuparrows}^e = \left( \begin{array}{cc}
0 & (Q_{\downuparrows}^e)_-   \\
(Q_{\downuparrows}^e)_+ & 0   \end{array} \right)
\end{eqnarray*}
and we have $(Q_{\downuparrows}^e)_+^* = (Q_{\downuparrows}^e)_- $ since $Q_{\downuparrows}^e$ are self-adjoint families in their domains. These local families transform under the action of the shift operators on $\mathbb{T}_s \times (\mathbb{T}_{\uparrow} \cap \mathbb{T}_{\downarrow} )^{(1)} \times M$ by 
\begin{eqnarray*}
	(1 \otimes S) Q^e_{\downarrow}(s,\phi,p) (1 \otimes S)^{-1} =  Q^e_{\uparrow}(s,\phi - 2 \pi,p).
\end{eqnarray*}
In particular, the local families $Q^e_{\downarrow}$ and $Q^e_{\uparrow}$  are mapped to each other over $\mathbb{T}_s \times (\mathbb{T}_{\uparrow} \cap \mathbb{T}_{\downarrow} )^{(1)} \times M$ under the adjoint action of the local Hilbert bundle isomorphisms $u_a$  acting on the component $\mathcal{F}$ in the fibres of the local Hilbert bundles, recall 1.1. \3

\noindent \textbf{Proposition 5.} The approximate sign $F^e_{\downuparrows} = Q^e_{\downuparrows} (1 + (Q^e_{\downuparrows})^2)^{-1/2}$ defines an element in the twisted K-theory group $K^0(\mathbb{T}^2 \times M, \tau)$. \3

\noindent Proof. According to the discussion above, the families $F^e_{\downuparrows}$ transform properly under the adjoint action of the local isomorphisms $u_a$ associated to the gerbe twisted by a decomposable class $\tau$. It remains to check the properties (A1-2-3) in the appendix. A straightforward computation using $\{\psi^1_k, \psi^0_l\}= 0 $ and $[e_k, f_l] = 0$ for all $k,l$, and the formula for the square of the supercharge $Q^e_{\downuparrows}$ in 1.6  gives 
\begin{eqnarray}\label{kernel}
	(Q^e_{\downuparrows}(s,\phi,p))^2 = l + l_0^e + l_0^f + (e_0 + \frac{\phi}{2 \pi})^2 + (f_0 + \frac{s}{2 \pi})^2,
\end{eqnarray} 
where $l^e_0$ and $l^f_0$ are given as in \eqref{l} and $l$ is the operator
\begin{eqnarray*}
l = \sum_{n \geq 1} n( \psi^0_{n} \psi^0_{-n} + \psi^1_{n} \psi^1_{-n}).
\end{eqnarray*}
The operator $l + l_0^e + l_0^f$ is a positive unbounded operator, which can be understood as the Laplacian of the supercharge. The eigenspaces of the local families $(Q^e_{\downuparrows}(s,\phi,p))^2$ are finite dimensional and the eigenvalues go to infinity. So the property A1 holds. It is also clear that A2 holds because $||\psi^0_0|| = ||\psi^1_0|| = 1$. Finally A3 holds because the finite particle states provide a dense domain for $Q^e$. \5 $\square$ \3

\noindent \textbf{3.4.} We can compute the index character of $Q^e_{\downuparrows}$ by taking the infinite time limit as in 2.3. The conventions for the superconnection analysis are as in 2.3. \3

\noindent \textbf{Theorem 4.} The index character of $Q^e$ associated to the vacuum twist $\xi$ is represented by 
\begin{eqnarray}\label{character34}
\underline{\textbf{ch-ind}}_0^{\tau}(\mathbb{F}) = \frac{ds}{2 \pi} \wedge \frac{d \phi}{2 \pi} \wedge \ch(F_{\xi}). 
\end{eqnarray}

\noindent \text{Proof.} We can perform another infinite time limit calculation to find a distribution valued form. For the computation we fix a one parameter family of superconnections $\mathbb{A}(t) = \sqrt{t} \pi^*(Q^e) + \pi^*(\nabla)$ over the covering space $\mathbb{T}_s \times \mathbb{R} \times M$. The computation itself follows the strategy used in 2.3. Let $\mathbb{F}(t) = t(\pi^*(Q^e)^2 + K_t)$ be the supercurvature so that 
\begin{eqnarray*}
K_t = \frac{1}{\sqrt{t}} \Big( \frac{d s}{2\pi} \psi_0^0 + \frac{d x}{2 \pi} \psi_0^1 \Big) + \frac{1}{t} \pi^*(F).
\end{eqnarray*}
The Volterra series gives 
\begin{eqnarray*}
	e^{-\mathbb{F}(t)} =  \sum_{n \geq 0} \frac{(-t)^n}{n!} K_t^n  e^{-t \pi^*(Q^e)^2}.
\end{eqnarray*}
Then we can use the binomial formula to simplify $K_t^n$ and after the contributions in the kernel of the supertrace are set to zero we are left with 
\begin{eqnarray*}
e^{-\mathbb{F}(t)} &\simeq& - t \sum_{n \geq 2} \frac{ds \wedge d x }{2 \pi^2} \wedge \frac{(- \pi^*(F))^{n-2}}{(n-2)!} \psi_0^0 \psi_0^1 e^{-t\pi^*(Q^e)^2} \\ 
&=& - \frac{t}{2 \pi^2} \psi_0^0 \psi_0^1  ds \wedge d x \wedge  e^{- \pi^*(F)} e^{-t\pi^*(Q^e)^2}. 
\end{eqnarray*}
Using \eqref{kernel}, we can write it as a product of the local families 
\begin{eqnarray*}
	&&\sqrt{\frac{t}{\pi}} e^{-t (e_0 + \frac{s}{2 \pi})^2}, \5 \sqrt{\frac{t}{\pi}}  e^{-t (f_0 + \frac{x}{2\pi})^2}  \5 \text{and} \\
	&& - \frac{1}{2 \pi} \psi_0^0 \psi_0^1 ds \wedge d x \wedge  e^{- \pi^*(F)} e^{-t (l + l_0^e + l_0^f)}
\end{eqnarray*}
each converges in the limit $t \rightarrow \infty$ resulting
\begin{eqnarray*}
	\lim_{t \rightarrow \infty} \text{sTr} (e^{-\mathbb{F}(t)})  =  \text{Tr}( P' \frac{i }{4 \pi} ds \wedge dx \wedge e^{-\pi^*(F)} \delta\Big(e_0 + \frac{s}{2\pi}\Big) \delta \Big(f_0 + \frac{x}{2 \pi}\Big)P')
\end{eqnarray*}
where $P'$ denotes the orthogonal projection onto the subspaces
\begin{eqnarray*}
\mathcal{H}^e_{\text{vac}} \otimes \bigoplus_{k \in \mathbb{Z}} S^k |0 \rangle \otimes \mathbb{C}^n
\end{eqnarray*}
in the fibres of the Fock bundle; $\mathcal{H}^e_{\text{vac}}$ denotes the 2-dimensional spinor vacuum which is the kernel of $l$. 

Now we can choose a section $\psi$ of $\pi : \mathbb{T} \times \mathbb{R} \times M \rightarrow \mathbb{T}^2 \times M$ and pull back the character and normalize it. Its equivalence class can be represented by \eqref{character34}. These steps are analogous to those in 2.3. \5 $\square$ \3

\noindent \textbf{3.5.} We can also compute the index in the twisted cohomology theory following the superconnection technique of \cite{BG11}. An even superconnection is a descent datum $\mathbb{A} = (\mathbb{A}_{\downuparrows}, \nabla_{\lambda}, \Omega_{\downuparrows})$ where $\mathbb{A}_{\downuparrows}$ are locally defined even superconnections, which are odd elements in
\begin{eqnarray*}
\mathbb{A}_{\downuparrows} \in \Lambda^{*}(\mathbb{T}_s \times \mathbb{T}_{\downuparrows} \times M, \text{End}(\textbf{S}^e \otimes (\textbf{F}^s \boxtimes \textbf{F}_{\downuparrows}) \otimes \xi ))
\end{eqnarray*}
 and satisfy the Leibnitz rule 
\begin{eqnarray*}
\mathbb{A}_{\downuparrows}(\Phi f) = (d  \Phi) \wedge f + (-1)^k \Phi \wedge \mathbb{A}_{\downuparrows} f
\end{eqnarray*}
for $\Phi \in \Lambda^k(\mathbb{T}_s \times \mathbb{T}_{\downuparrows} \times M)$ and $f \in \Lambda^{*}(\mathbb{T}_s \times \mathbb{T}_{\downuparrows} \times M, \textbf{S}^e \otimes (\textbf{F}^s \boxtimes \textbf{F}_{\downuparrows}) \otimes \xi )$ and are subject to the transformations
\begin{eqnarray*}
u_a^*(\mathbb{A}_{\uparrow}) = \mathbb{A}_{\downarrow} \otimes \text{id} + \text{id} \otimes \nabla_{\lambda} \5 \text{on} \5 \mathbb{T}_s \times ( \mathbb{T}_{\downarrow} \cap \mathbb{T}_{\uparrow})^{(1)} \times U_a. 
\end{eqnarray*}
The $\mathbb{Z}_2$ grading  (to the odd and even elements) is given by the total degree of the $\mathbb{Z}_2$ grading in the complex of differential forms and the $\mathbb{Z}_2$ grading of the endomorphism bundle coming from the polarization of $\textbf{S}^e$. 

By slightly abusing the notation, we will denote by $\Omega_{\downuparrows}$ the forms introduced in 1.1 extended on $\mathbb{T}_s \times \mathbb{T}_{\downuparrows} \times M$. Then we use the local supercurvatures associated to $\mathbb{A}_{\downuparrows}$ to define the global curvature form
\begin{eqnarray*}
\mathbb{F}_{\downuparrows} = \mathbb{A}_{\downuparrows}^2 + \Omega_{\downuparrows}
\end{eqnarray*}
and define its character by 
\begin{eqnarray*}
\textbf{ch-ind}_0^{\tau}(\mathbb{F}) = \varphi(\text{sTr}(e^{- \mathbb{A}_{\downuparrows}^2 - \Omega_{\downuparrows}})).
\end{eqnarray*}
Now we apply the supertrace in the even formalism. This defines a cocycle in the twisted cohomology determined by the twisted differential $d + H$, and its twisted cohomology class is independent on the choice of the superconnection, see \cite{BG11}. 

The following can be easily adapted using the idea in the corresponding result of 1.6 in the odd twisted K-theory, and the computation in 3.4. \3

\noindent \textbf{Proposition 6.} The character $\textbf{ch-ind}_0^{\tau}(\mathbb{F})$ is represented by 
\begin{eqnarray*}
\frac{ds}{2 \pi} \wedge \frac{d\phi}{2 \pi} \wedge \ch(F_{\xi})
\end{eqnarray*}
in the twisted cohomology group $H^{\text{even}}(\mathbb{T}^2 \times M, H)$.

\section{The Case of a Torus}

\textbf{4.1.} We use the torus $\mathbb{T}_{\phi} \times \mathbb{T}^n$ to demonstrate how the supercharge techniques can be applied. Fix the generators $\{d \theta_i: 1 \leq i \leq n\}$ for the cohomology ring of $\mathbb{T}^n$. To simplify the notation, we suppose that these forms are normalized so that they integrate to $1$ over the circles: $\int_{\mathbb{T}_i} d \theta_i = 1$. Without a loss of generality we can fix a twisting line bundle $\lambda$ with a Chern class equal to $k d \theta_1 \wedge d \theta_2$ where $k$ can be any nonzero integer. Let $\Lambda_n$ denote the exterior algebra over $\mathbb{Z}$ with $n$ generators with $\Lambda_0 = \Lambda_0^{\text{even}} = \mathbb{Z}1$. \3

\noindent \textbf{Proposition 7.} Let $\tau$ denote a cup product 3-cohomology class on a torus $\mathbb{T}_{\phi} \times \mathbb{T}^n$ with $n \geq 2$. The twisted K-theory groups are isomorphic to
\begin{eqnarray*}
K^0(\mathbb{T}^{n+1}, \tau) \simeq \Lambda_{n-2}^{\text{odd}} \oplus \Lambda^{\text{odd}}_{n-1} \oplus \frac{\Lambda_{n}^{\text{odd}}}{k\Lambda_{n-2}^{\text{odd}}}, \5 K^1(\mathbb{T}^{n+1}, \tau) \simeq \Lambda_{n-2}^{\text{even}} \oplus \Lambda^{\text{even}}_{n-1} \oplus \frac{\Lambda_{n}^{\text{even}}}{k\Lambda_{n-2}^{\text{even}}}. 
\end{eqnarray*}

\noindent Proof. Recall that we have the isomorphism \eqref{twistedk}. 

Let us start by computing all $x \in K^{*}(\mathbb{T}^n)$ which satisfy $(1 - \lambda) \otimes x = 0$. The K-theory ring on a torus is isomorphic to the exterior algebra over $\mathbb{Z}$. Identify the generators of the exterior algebra with the forms $d\theta_1, \ldots, d \theta_n$. The element $\lambda$ is represented by $1 - k d \theta_1 \wedge d \theta_2$. Then it is sufficient to find the forms $\Omega$ which verify 
\begin{eqnarray}\label{invariance}
 k d \theta_1 \wedge d \theta_2 \wedge \Omega = 0. 
\end{eqnarray}
Clearly we just need to find all the wedges which contain one of the forms $d \theta_1$ and $d \theta_2$. 

The invariant 1-forms (i.e. those which verify \eqref{invariance}) are integer multiplets of $d \theta_1$ and $d \theta_2$. An invariant k-form is a sum of forms of  the following 2 types: $d \theta_1 \wedge \Theta$ where $\Theta$ is any $(k-1)$-form in the algebra over $\mathbb{Z}$ generated by $\{d \theta_2, \ldots, d \theta_n\}$, and $d \theta_2 \wedge \Theta'$ where $\Theta'$ is any $(k-1)$-form in the algebra over $\mathbb{Z}$ generated by $\{d \theta_3, \ldots, d \theta_n\}$. Therefore, the invariant subspace in the odd integer cohomology is isomorphic to $\Lambda^{\text{even}}_{n-1} \oplus \Lambda^{\text{even}}_{n-2}$. Similarly, the invariant subspace in the even cohomology is isomorphic to $\Lambda^{\text{odd}}_{n-1} \oplus \Lambda^{\text{odd}}_{n-2}$. The invariant parts in the twisted K-theory groups are isomorphic to
\begin{eqnarray*}
\{x \in K^1(\mathbb{T}^n) : (1 - \lambda) \otimes x = 0 \} &=& \Lambda^{\text{even}}_{n-1} \oplus \Lambda^{\text{even}}_{n-2} \5 \text{and} \\
\{x \in K^0(\mathbb{T}^n) : (1 - \lambda) \otimes x = 0 \} &=& \Lambda^{\text{odd}}_{n-1} \oplus \Lambda^{\text{odd}}_{n-2}.
\end{eqnarray*}
Both groups are torsion free, and so the group extension of \eqref{twistedk} is trivial.

Similarly, by studying an exterior algebra one can check that the normal subgroup $(1 - \lambda) \otimes K^0(M)$ is isomorphic to $k \Lambda^{\text{even}}_{n-2}$ and $(1 - \lambda) \otimes K^{1}(\mathbb{T}^n)$ is isomorphic to $k \Lambda^{\text{odd}}_{n-2}$. \5 $\square$ \3

\noindent \textbf{4.2.} We first study the odd twisted K-theory of $\mathbb{T}^{n+1}$ with $n \geq 2$. The supercharge construction of \cite{HM12}, which was sketched in 1, provides a concrete realization for all the twisted K-theory classes associated with the subgroup isomorphic to
\begin{eqnarray}\label{torsion}
\frac{K^0(\mathbb{T}^n)}{(1 - \lambda) \otimes K^0(\mathbb{T}^n)} \leq K^1(\mathbb{T}^{n+1}, \tau). 
\end{eqnarray}
Recall that we have the index character for this class given in 1.4. We can modify the supercharge of 1.3 slightly by creating a negative spectral flow around $\mathbb{T}_{\phi}$ as well. Change the local Fock bundle isomorphism by $u_a^{-1}(\phi,p) = \overline{s_a}(p) S^{-1}$ (recall \eqref{locuni}) which leads to the following modification to the supercharge
\begin{eqnarray*}
Q_{\downuparrows} = \sum_k \psi_k \otimes e_{-k} - \frac{\phi}{2 \pi} \psi_0 \otimes \textbf{1}.
\end{eqnarray*}
We also need to change the twisting line bundle $\lambda$ to its complex conjugate bundle. The Dixmier-Douady cohomology class of the gerbe remains invariant because the negative spectral flow and the complex conjugate twisting bundle send the associated cohomology classes to their additive inverses. Now the index character can be represented by  
\begin{eqnarray*}
\underline{\textbf{ch-ind}}^{\tau}_1(\mathbb{F}) = - \sqrt{\pi} \frac{d \phi}{2 \pi} \wedge \ch(F_{\xi}). 
\end{eqnarray*}
Since the Chern character $\ch$ is an injective group homomorphism in the case of ordinary K-theory, we find all the classes of type \eqref{torsion} by tuning the vacuum bundle $\xi$ and by using the group structure.\3 

\noindent \textbf{4.3.} The suspension to the even twisted K-theory can be used to construct analytic realizations for the elements in the subgroup 
\begin{eqnarray}\label{torsion2}
\frac{K^1(\mathbb{T}^n)}{(1 - \lambda) \otimes K^1(\mathbb{T}^n)} \leq K^0(\mathbb{T}^{n+1}, \tau). 
\end{eqnarray}
We can follow the construction of the section 3 and define a supercharge with spectral flow around the circles $\phi$ and $\theta_i$ for any $1 \leq i \leq n$. Again we can allow both positive and negative spectral flow in the direction $\phi$  which gives the following supercharges 
\begin{eqnarray*}
	Q_{\downuparrows}^e &:& \mathbb{T}_{\downuparrows} \times \mathbb{T}^n \rightarrow \textbf{Fred}^{(0)}_{\Psi} \\
	Q_{\downuparrows}^e(\phi,\theta_1, \ldots, \theta_n) &=& \sum_k \psi^0_k \otimes e_{-k} + \sum_k \psi^1_k \otimes f_{-k} \pm \phi \psi^0_0 \otimes \textbf{1} + \theta_i \psi^1_0  \otimes \textbf{1}. 
\end{eqnarray*}
Suppose that $\xi$ is a complex vector bundle over $\mathbb{T}_{\theta_1} \times \cdots \mathbb{T}_{\theta_{i-1}} \times \mathbb{T}_{\theta_{i+1}} \times \cdots \times \mathbb{T}_{\theta_{n}}$ which we use to tensor the vacuum of the local Hilbert bundles. The characters of the supercharges are now
\begin{eqnarray*}
\underline{\textbf{ch-ind}}^{\tau}_0(\mathbb{F}) = \pm \frac{d \phi}{2 \pi} \wedge d \theta_i \wedge \ch(F_{\xi}).
\end{eqnarray*}
Again, by tuning the vacuum bundle $\xi$ we find representatives for all the elements in the subgroup of twisted K-theory that is isomorphic to \eqref{torsion2}. In particular, we find all the torsion elements in the twisted K-theory on tori.

\section*{Appendix} 

\textbf{A.1 - twisted K-theory.} We apply the twisted K-theory model introduced by Atiyah and Segal in \cite{AS04}. Let $X$ denote a compact manifold with a good open cover $\{U_i\}$. The twist class $\tau$ is fixed by choosing a local family of Hilbert bundles $\textbf{H}_i \rightarrow U_i$ together with local isomorphisms $\varphi_{ij}: \textbf{H}_i \otimes \lambda_{ij} \rightarrow \textbf{H}_j$ where $\lambda_{ij} \rightarrow U_{ij}$ are complex line bundles. The Dixmier-Douady class of the gerbe is the Cech-cohomology class in $H^2(X, \underline{\mathbb{T}})$ whose components are determined by the local isomorphisms 
\begin{eqnarray*}
f_{ijk}: \lambda_{ij} \otimes \lambda_{jk} \rightarrow \lambda_{ik} 
\end{eqnarray*}
over the 3-intersections $U_{ijk}$. The 3-cohomology class $\tau$ of the twist is the image of $f$ under the group isomorphism $H^2(X, \underline{\mathbb{T}}) \rightarrow H^3(X, \mathbb{Z})$ associated to the group extension $\mathbb{Z} \rightarrow \mathbb{R} \rightarrow \mathbb{T}$. 

Let $\mathcal{B}(\hil)$ denote the space of bounded operators on an infinite dimensional separable complex Hilbert space $\hil$. We equip   $\mathcal{B}(\hil)$ with the strong operator topology. Then $\textbf{Fred}^{(0)}$ denotes the space of operators $A \in \mathcal{B}(\hil)$ such that there is $B \in \mathcal{B}(\hil)$, a parametrix of $A$, so that $1 - AB$ and $1 -BA$  are compact operators. The space of compact operators is equipped with the norm topology. Denote by $\textbf{Fred}^{(1)}$ the subspace of self-adjoint elements with both positive and negative essential spectrum. 

Associated to the local Hilbert bundles $\textbf{H}_i$ there are bundles of $\textbf{Fred}^{\bullet}$ operators for $\bullet = 0,1$. Under the conjugation action by the local unitary families $\varphi_{ij}$, these local bundles glue together to form a Fredholm bundle over $X$. By definition, the elements of the twisted K-theory groups $K^{\bullet}(X, \tau)$ are homotopy classes of continuous sections of the Fredholm bundles for $\bullet = 0,1$. Precisely speaking, twisted K-theory is twisted by a bundle-realization of the cohomology class $\tau$ and there is no canonical way to identify two different twisted K-theory groups associated with two isomorphic, but different realizations for $\tau$. However, this should not lead to any confusion since the gerbe is always carefully defined when twisted K-theory is discussed.  

We will work with sections of a bundle of unbounded Fredholm operators. The spaces of unbounded Fredholm operators $\textbf{Fred}^{\bullet}_{\Psi}$ are equipped with the topology induced by the map 
\begin{eqnarray*}
\Psi: Q \mapsto \frac{Q}{\sqrt{1 + Q^2}} = F_Q. 
\end{eqnarray*}
The image $F_Q$ is called the approximate sign of $Q$. In the case of even twisted K-theory, we work with $\mathbb{Z}_2$-graded Hilbert bundles and with antidiagonal self-adjoint Fredholm sections $Q_i: \textbf{H}^{\pm}_i \rightarrow \textbf{H}^{\mp}_i$. Then the component $(F_Q)^+_i: \textbf{H}^{+}_i \rightarrow \textbf{H}^{-}_i$ will realize the class in $K^0(X, \tau)$. Given an unbounded Fredholm section $Q$ it is proved in \cite{FHT13} that the approximate sign defines an element in twisted K-theory if the following conditions are satisfied:  
\begin{quote}
\textbf{A1.} $Q$ has eigenvalues of finite multiplicity and they go to infinity. 

\textbf{A2.} Each $x \in X$ has a neighborhood $U$ such that $Q(y) = Q' + b(y)$ if $y \in U$ where $Q'$ is the constant family $Q'(y) = Q(x)$ for all $y \in U$ and $b$ is a family of bounded operators over $U$ which is continuous in the norm topology. 

\textbf{A3.} The fibre $\hil$ at $x \in X$ has a dense subspace $H \subset \hil$ so that $Q(x)$ is densely defined for all $x \in X$. 
\end{quote}
The property A1 is sufficient to prove that $1 - (F_Q)^2(x)$ is a compact operator and therefore $F_Q$ is its own parametrix. According to the proof of Corollary 3.39 in \cite{FHT13}, A2 implies that $x \mapsto 1 - (F_Q)^2(x)$ is continuous in the norm topology. If A1 and A2 hold, then A3 is a sufficient condition for the continuity of $x \mapsto F_Q(x)$ in the strong operator topology (see Lemma 3.43 in\cite{FHT13}). \3

\noindent \textbf{A.2 - twisted cohomology.} Suppose that a gerbe and so a representative of the twisting class is fixed. Let $H$ denote the 3-curvature, or an $H$-flux, of the gerbe. Then the twisted cohomology is the cohomology of the cocycle complex $(\Lambda^{*}(X), d + H)$. If we have two realizations for the twisting class with the 3-curvatures $H_1$ and $H_2$, then there is an isomorphism between the cohomology groups determined by the differentials $d+H_1$ and $d+H_2$ which is given by 
\begin{eqnarray*}
\Omega \mapsto e^{- \Phi} \wedge \Omega,
\end{eqnarray*}
if $\Phi$ is a form satisfying $H_2 = H_1 + d \Phi$. 

\bibliographystyle{plain}

\end{document}